\documentclass[11pt,a4paper]{amsart}
\usepackage{amssymb,amsmath}
\usepackage[american,french]{babel}
\usepackage{amsopn}
\usepackage{graphicx}
\usepackage{fancyhdr}
\usepackage[T1]{fontenc}
\usepackage{mathrsfs}
\usepackage{ulem} 
\usepackage{color}

\newcommand{\rr}{\mathbb{R}}
\newcommand{\eps}{\varepsilon}
\newcommand{\nn}{\mathbb{N}}
\newcommand{\cc}{\mathbb{C}}

\numberwithin{equation}{section}

\begin{document}
\newcounter{equa}
\selectlanguage{american}

\newtheorem{theorem}{Theorem}[section]
\newtheorem{corollary}[theorem]{Corollary}
\newtheorem{lemma}[theorem]{Lemma}
\newtheorem{definition}[theorem]{Definition}
\newtheorem{proposition}[theorem]{Proposition}
\newtheorem{remark}[theorem]{Remark}

\title[Semiclassical hypoelliptic estimates]{Semiclassical hypoelliptic estimates with a big loss of derivatives}

\author{Alberto Parmeggiani, Karel Pravda-Starov}

\address{\noindent \textsc{Dipartimento di Matematica,
Universit\`a di Bologna,
Piazza di Porta S. Donato, 540126 Bologna, Italia}}
\email{alberto.parmeggiani@unibo.it}

\address{\noindent \textsc{
Universit\'e de Cergy-Pontoise,
CNRS UMR 8088,
D\'epartement de Math\'ematiques,
95000 Cergy-Pontoise, France}}
\email{karel.pravda-starov@u-cergy.fr}


\begin{abstract}
We study the pseudospectral properties of general pseudodifferential operators around a doubly characteristic point and provide necessary and sufficient conditions for semiclassical hypoelliptic a priori estimates with a big loss of derivatives to hold.
\end{abstract}

\maketitle

\section{Introduction}
Over the past years, there has been a renewed interest in the analysis of the spectra and resolvents of non-selfadjoint operators with double characteristics. 
This interest finds some of its grounds in the study of the long-time behavior of evolution equations associated with non-selfadjoint operators 
$$\left\lbrace
\begin{array}{c}
(\partial_t+P) u(t,x)=0  \\
u(t,\textrm{\textperiodcentered})|_{t=0}=u_0.
\end{array} \right.$$
This is for instance the case in the analysis of kinetic equations and the study of the trend to equilibrium in statistical physics.

The study of doubly characteristic operators has a long and distinguished tradition in the analysis of partial differential equations \cite{H,hormander,S}. 
The simplest examples of such operators are given by quadratic differential operators
$$Q(x,D_x)=\sum_{|\alpha+\beta|=2}q_{\alpha,\beta}x^{\alpha}D_x^{\beta}, \quad x \in \rr^n,$$
with $q_{\alpha,\beta} \in \cc$, $D_{x_j}=i^{-1}\partial_{x_j}$, $\alpha, \beta \in \nn^n$. In the elliptic case, the spectrum of these operators has been understood and 
described explicitly for some time~\cite{S}. On the other hand, the pseudospectral study of these operators is much more recent. Studying the pseudospectrum of an operator 
is studying the level lines of the norm of its resolvent 
$$\mathrm{Spec}_{\varepsilon}(A)=\Big\{z \in \cc;  \ \|(A-z)^{-1}\| \geq \frac{1}{\varepsilon} \Big\}, \quad \eps>0,$$
with the convention $\|(A-z)^{-1}\|=+\infty$ if $z$ belongs to the spectrum~$\mathrm{Spec}(A)$ of $A$. 
The breadth of these sets allows to analyze the spectral stability of the operator under small perturbations. Indeed, the pseudospectrum may be defined in an equivalent 
way~\cite{roch} in terms of the spectra of the operator perturbations 
$$\mathrm{Spec}_{\eps}(A)=\bigcup_{B \in \mathscr{L}(H), \ \|B\|_{\mathscr{L}(H)} \leq \eps}{\mathrm{Spec}(A+B)},$$
where $\mathscr{L}(H)$ stands for the set of bounded linear operators on $H$. 
The pseudospectral study of a variety of operators has received much recent interest in a diverse array of problems. For further details and motivations, 
we refer the reader to the overview of this topic presented in the book~\cite{trefethen2}, and to all the references therein. For now, let us simply notice that 
the study of the pseudospectrum is non-trivial only for non-selfadjoint operators, or more precisely for non-normal operators. In fact, the classical formula 
\begin{equation}\label{1}
\forall z \not\in \mathrm{Spec}(A), \quad \|(A-z)^{-1}\| = \frac{1}{\textrm{dist}(z,\mathrm{Spec}(A))}, 
\end{equation}
emphasizes that the resolvent of a normal operator cannot blow up far from its spectrum, and that the spectrum is stable under small perturbations 
\begin{equation}\label{2}
\mathrm{Spec}_{\eps}(A) = \{z \in \cc ;\  \textrm{dist}(z,\mathrm{Spec}(A)) \leq \eps\}. 
\end{equation}
However, formula (\ref{1}) does not hold anymore for non-normal operators and the behavior of the resolvent for such operators can be intricate by becoming very 
large in norm far from the spectrum. As a consequence, the spectrum of these operators may be very unstable under small perturbations. The rotated harmonic oscillator
$$P=-\partial_x^2+e^{i \theta} x^2, \quad -\pi <\theta <\pi, \ \theta \neq 0,$$
is a noticeable example of elliptic quadratic operator whose spectrum is very unstable under small perturbations. The seminal works~\cite{boulton,davies}
have indeed shown that its resolvent $\|(P-z)^{-1}\|$ exhibits a rapid growth in some regions of the resolvent set far away from the spectrum, and that some strong 
spectral instabilities are developing in some regions with a specific geometry, which have been sharply described in the works~\cite{boulton,kps13}. These phenomena of 
spectral instabilities are not peculiar to the rotated harmonic oscillator. They were shown to be the typical behavior of any non-normal elliptic quadratic 
operator~\cite{kps14,kps17,kps19}, with a rapid resolvent growth along any ray lying inside the range of the Weyl symbol of these operators.
This is linked to some properties of microlocal non-solvability and to violations of the adjoint condition to the so-called Nirenberg-Treves condition $(\Psi)$,
which allow the construction of quasimodes~\cite{daviessemi,dencker,hormander,kps12,kps14,zworski1,zworski2}. Similar types of spectral instabilities were shown to 
occur for general pseudodifferential operators around a doubly characteristic point, when the quadratic approximations of these operators at the doubly characteristic 
set are non-normal~\cite{kps15}. Starting from these early insights, there has been a series of recent works~\cite{dencker,HeHiSj2,HeHiSj3,HSS06,kps3,kps4,kps15,sjosemi,viola1,viola2} 
intending to provide a sharp description of the spectral and pseudospectral properties of general pseudodifferential operators around a doubly characteristic point. In the present work, 
we aim at completing this picture and at describing the pseudospectral behavior of a general pseudodifferential operator around a doubly characteristic point by refining the understanding of the 
underlying geometry ruling these phenomena.

\section{Setting of the analysis}\label{setting}

Let $m(\cdot;h) : \rr^{2n} \longrightarrow ]0,+\infty[$ be an order function (see Dimassi-Sj\"ostrand's book \cite{DS}), that is, 
$$\exists C_0,N_0>0,\forall 0<h \leq 1, \forall X,Y\in \rr^{2n}, \quad m(X;h)\leq C_0 \langle{X-Y\rangle}^{N_0} m(Y;h),$$
with $\langle X \rangle=(1+|X|^2)^{1/2}$, where $|\cdot|$ is the Euclidean norm. We consider the symbol class of $h$-dependent symbols whose growth
is controlled by the order function $m$ given by
\begin{multline}\label{sl70}
S(m)=\{ a(\cdot;h)\in C^{\infty}(\rr^{2n},\cc); \forall \alpha \in \nn^{2n}, \exists C_{\alpha}>0, \\ \forall  0<h \leq 1, \forall X \in \rr^{2n},
\quad  |\partial_X^{\alpha} a(X;h)| \leq C_{\alpha}m(X;h)\}.
\end{multline} 
In the present work, we study a semiclassical pseudodifferential operator
\begin{equation}\label{sl4}
P=p^w(x,hD_x;h)=\frac{1}{(2\pi)^n}\int_{{\rr}^{2n}}e^{i(x-y)\cdot\xi}p\Big(\frac{x+y}{2},h\xi;h\Big)u(y)dyd\xi,
\end{equation}
defined by the semiclassical Weyl quantization of a symbol $p(x,\xi;h)$ admitting a semiclassical asymptotic expansion in the symbol class~$S(1)$,
\begin{equation}\label{sl5}
p(x,\xi;h) \sim \sum_{j=0}^{+\infty}  p_j(x,\xi)h^j.
\end{equation}
The symbols $p_j \in S(1)$ in the asymptotic expansion are supposed to be independent of  
the semiclassical parameter $0<h \leq 1$.
We assume that the real part of the principal symbol is nonnegative 
\begin{equation}\label{sl6}
\textrm{Re }p_0(X) \geq 0, \quad X=(x,\xi) \in \rr^{2n},
\end{equation}
and elliptic at infinity 
\begin{equation}\label{sl7}
\exists C>1, \forall |X|\geq C, \quad  \textrm{Re }p_0(X) \geq \frac{1}{C}.
\end{equation}
These two assumptions imply that there exists a neighborhood of zero in the complex plane such that the analytic family of bounded operators 
$$P-z: L^2(\rr^n)\longrightarrow L^2(\rr^n), \quad z\in {\rm neigh}(0,\cc),$$
is Fredholm of index $0$, when the semiclassical parameter  $0<h\ll 1$ is small enough~\cite{dencker}. An application of the analytic Fredholm theory 
shows that the spectrum of the operator $P$ in a small neighborhood $V$ of $0$, that we may take of the form $V=D(0,c)$ (the open disk in $\mathbb{C}$ 
centered at $0$ of radius $c$), with $0<c \leq 1$, is discrete and is uniquely composed of eigenvalues with finite algebraic multiplicity.

We assume further that the characteristic set of the real part of the principal symbol is reduced to a single point
\begin{equation}\label{sl8}
(\textrm{Re }p_0)^{-1}(\{0\})=\{0\} \subset \rr^{2n}
\end{equation}
and that this point is doubly characteristic for the principal symbol $p_0$
\begin{equation}\label{sl9}
p_0(0)=\nabla p_0(0)=0,
\end{equation}
so that we may write
\begin{equation}\label{sl10}
p_0(Y)=q(Y)+\mathcal{O}(Y^3),\quad Y\to 0,
\end{equation}
$q$ being the quadratic term in the Taylor expansion of the principal symbol at~$0$. 

We aim at studying the spectral and pseudospectral properties of the operator $P$ in a neighborhood of $0$. 
As mentioned above, the study of this problem was started in~\cite{kps3,kps4}, where the first lines of this spectral and pseudospectral picture were sketched out. 

The results of~\cite{kps3} actually provide a first localization of the spectrum of the operator $P$ in any given $h$-ball 
centered at $z=0$.  
More specifically, when the quadratic approximation of the principal symbol is elliptic on a particular vector subspace $S$ of 
phase space defined as its singular 
space\footnote{We refer the reader to the Appendix, Section~\ref{appendix}, for miscellaneous facts about quadratic operators and the definition of the singular space.}
\begin{equation}\label{kps223}
X \in S, \quad q(X)=0 \Longrightarrow X=0,
\end{equation}
then, for any given constant $C>1$ and any fixed neighborhood $\Omega \subset \cc$ of the spectrum $\mathrm{Spec}(q^w(x,D_x))$
of the quadratic operator $q^{w}(x,D_x)$ described in the Appendix (Section~\ref{appendix}),
there exist positive constants $0 < h_0 \leq 1$, $C_0>0$ such that for all $0<h \leq h_0$, $|z|\leq C$ satisfying
$$z-p_1(0) \notin \Omega,$$
we have
\begin{equation}\label{eq1.7}
h\|u\|_{L^2}\leq C_0\|(P-hz)u\|_{L^2}, \quad u \in \mathscr{S}(\rr^n),
\end{equation}
where $p_1(0)$ stands for the value of the subprincipal symbol at the doubly characteristic point $0 \in \rr^{2n}$. This result indicates that 
the spectrum of $P$ in any $h$-ball 
centered at $z=0$,
is localized in an $h$-neighborhood of the spectrum of its quadratic approximation shifted by the value of the 
subprincipal symbol at the doubly characteristic point
$$p_1(0)+\mathrm{Spec}(q^w(x,D_x)).$$
Under the same assumptions, this pseudospectral picture was completed by the following result about the spectrum~\cite{kps4}:  
For any given
$C>0$, there exists $0<h_0 \leq 1$, such that for all $0< h \leq h_0$, the spectrum of the operator $P$ in the open disk $D(0,Ch)$ is given by eigenvalues $z_k$ satisfying a semiclassical expansion of the type
\begin{equation}\label{eq1.9}
z_{k} \sim h (\lambda_{k} + p_1(0) + h^{1/N_{k}} \lambda_{k,1} + h^{2/N_{k}} \lambda_{k,2} +\ldots),
\end{equation}
where $\lambda_{k}$ is an eigenvalue of the quadratic operator $q^w(x,D_x)$ located in the fixed ball $D(0,C)$, $N_{k}$ is the dimension of the corresponding generalized eigenspace, and the
$\lambda_{k,j} \in \cc$ are some complex constants.

We next consider the remainder term in the principal symbol
\begin{equation}\label{hel1}
r(X)=p_0(X)-q(X), 
\end{equation}
and assume further the existence of a closed angular sector $\Gamma$ with vertex at 0, and a neighborhood $V$ of the origin in $\rr^{2n}$
such that 
\begin{equation}\label{re1}
r(V) \setminus \{0\} \subset \Gamma \setminus \{0\} \subset \{z \in \cc; \textrm{Re }z > 0\}.
\end{equation}
When the quadratic approximation $q^w(x,D_x)$ enjoys some subelliptic properties, sharp resolvent estimates may be derived outside an $h$-ball 
centered at $z=0$,
in a parabolic region with a particular geometry.
More specifically, when the quadratic form $q$ has a zero singular space, i.e. $S=\{0\}$, we consider the smallest integer $0 \leq k_0 \leq 2n-1$ satisfying
\begin{equation}\label{ning1}
\Big(\bigcap_{j=0}^{k_0}\textrm{Ker}\bigl(\textrm{Re }F(\textrm{Im }F)^j\bigr)\Big) \cap \rr^{2n}=\{0\},
\end{equation}
where $F$ is the Hamilton map of $q$ (see the Appendix, Section~\ref{appendix}).
It was shown in~\cite{kps4} that for any given sufficiently small constant $c_0 > 0$ there exist positive constants $0<h_0 \leq 1$, $C \geq 1$, $C_0>0$, such that for all
$0<h \leq h_0$, $u \in \mathscr{S}(\rr^n)$, and $z \in \Omega_{h}$,
\begin{equation}\label{ning4}
 h^{\frac{2k_0}{2k_0+1}}|z|^{\frac{1}{2k_0+1}}\|u\|_{L^2} \leq C_0 \|Pu-zu\|_{L^2},
\end{equation}
where
\begin{equation}\label{sl1}
\Omega_{h}=\Big\{z \in \cc;\ \textrm{Re }z \leq \frac{1}{C} h^{\frac{2k_0}{2k_0+1}}|z|^{\frac{1}{2k_0+1}}, \ Ch \leq |z| \leq c_0  \Big\}.
\end{equation}
The term $h^{\frac{2k_0}{2k_0+1}}|z|^{\frac{1}{2k_0+1}}$ increases when the spectral parameter $z$ moves away from the origin in the region where $Ch \leq |z| \leq c_0$.
When the spectral parameter is of magnitude $h$, we recover the semiclassical estimate (\ref{eq1.7}), and we emphazise that the resolvent estimate 
$$(P-z)^{-1} =\mathcal{O}(h^{-\frac{2k_0}{2k_0+1}}|z|^{-\frac{1}{2k_0+1}}): L^2(\rr^n)\longrightarrow L^2(\rr^n),$$
and the geometry of the parabolic region where it holds, are directly related to the loss of $2k_0/(2k_0+1)$ derivatives appearing in the global subelliptic estimate satisfied by the quadratic 
approximation of the operator at the doubly characteristic point (see the Appendix, Section~\ref{appendix}),
$$\|\langle(x,D_x)\rangle^{2/(2k_0+1)}u\|_{L^2} \leq C(\|q^w(x,D_x)u\|_{L^2}+\|u\|_{L^2}).$$
\begin{figure}[hhh]\label{des3.1}
\caption{The estimate $h^{\frac{2k_0}{2k_0+1}}|z|^{\frac{1}{2k_0+1}}\|u\|_{L^2} \leq C_0\|Pu-zu\|_{L^2}$ is fulfilled when $z$ belongs to the dark grey region of the figure, whereas the estimate 
$h\|u\|_{L^2} \leq C_0\|Pu-zu\|_{L^2}$ is fulfilled in the light grey one.}
\includegraphics[scale=0.6]{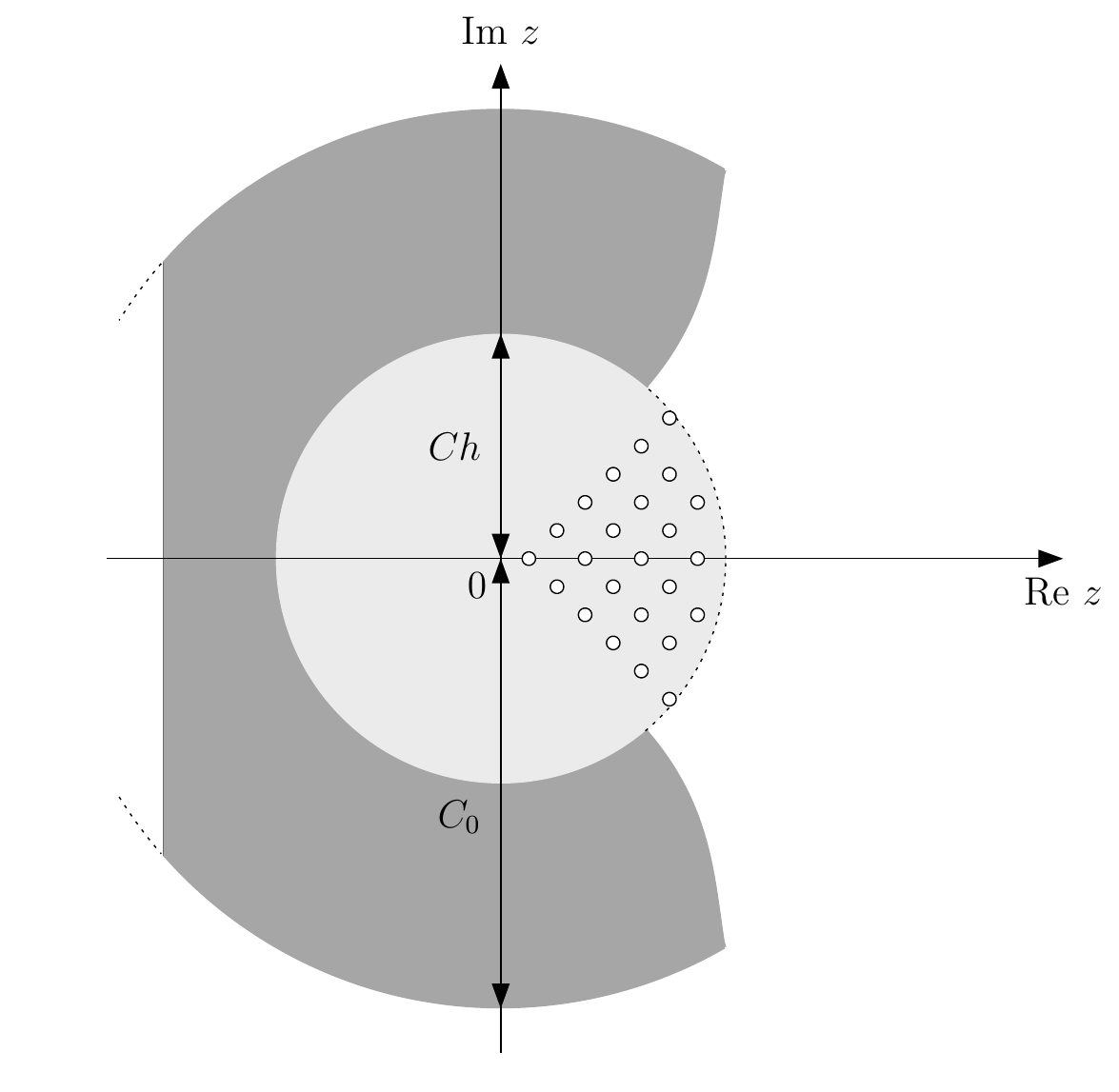}
\end{figure}
\hspace{-0.2cm}
These results show that the algebraic structure of the singular space allows to sharply account for the spectral and pseudospectral properties of pseudodifferential operators around a 
doubly characteristic point.
The picture drawn so far has been recently completed by Viola~\cite{viola1,viola2}. 
In these works, the author studies the case when the spectral parameter $z$ enters deeper into the numerical range and may grow slightly more rapidly than the semiclassical parameter~$h$ 
outside the parabolic region $\Omega_h$. His result shows that polynomial resolvent bounds still hold in a larger $h(\log \log h^{-1})^{1/n}$-neighborhood 
of $z=0$.
More precisely, under the previous assumptions with a zero singular space, Viola shows that for any given $\rho>0$, there exist positive constants $C_0,C_1>0$ such that the resolvent 
$$(P-z)^{-1} : L^2(\rr^n) \longrightarrow L^2(\rr^n),$$
exists and satisfies the bound
$$\|(P-z)^{-1}\|_{\mathscr{L}(L^2)}=\mathcal{O}(h^{-1-\rho}),$$
when $0<h \ll 1$, as long as the spectral parameter $z$ obeys
$$|z| \leq \frac{1}{C_0}h\Big(\log \log \frac{1}{h}\Big)^{1/n}, \quad \textrm{dist}(z,\mathrm{Spec}(q^w(x,D_x))) \geq h e^{-\frac{1}{C_1}(\log \log \frac{1}{h})^{1/n}}.$$ 
The next figure\footnote{Courtesy of Joe Viola.} is an illustration of a typical region in the complex plane where this resolvent estimate holds, for decreasing value of $h$.

\begin{figure}[ht]
\begin{minipage}[b]{0.45\linewidth}
\centering
\includegraphics[scale=0.4]{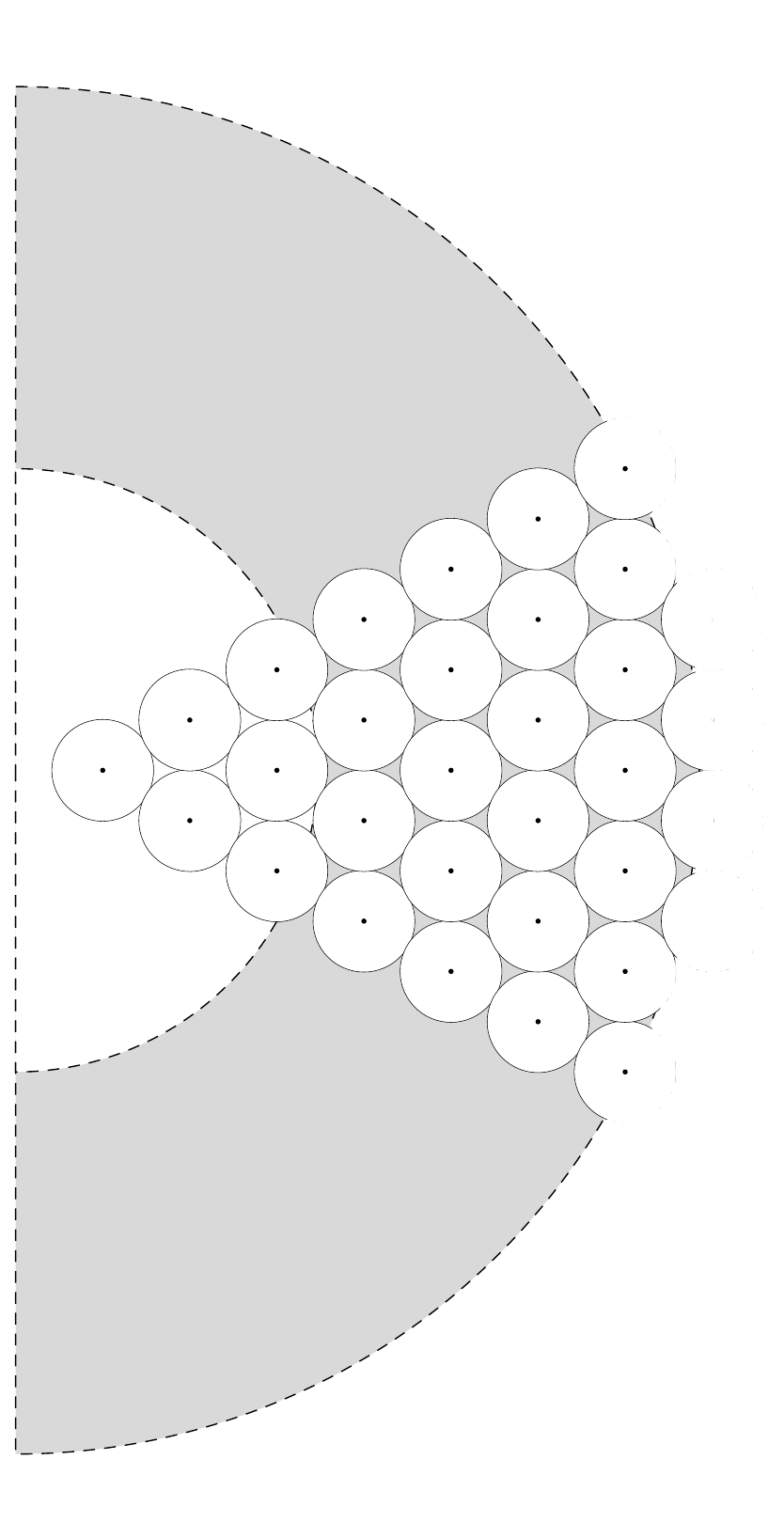}
\label{fig:figure1.1}
\end{minipage}
\hspace{0.5cm}
\begin{minipage}[b]{0.45\linewidth}
\centering
\includegraphics[scale=0.4]{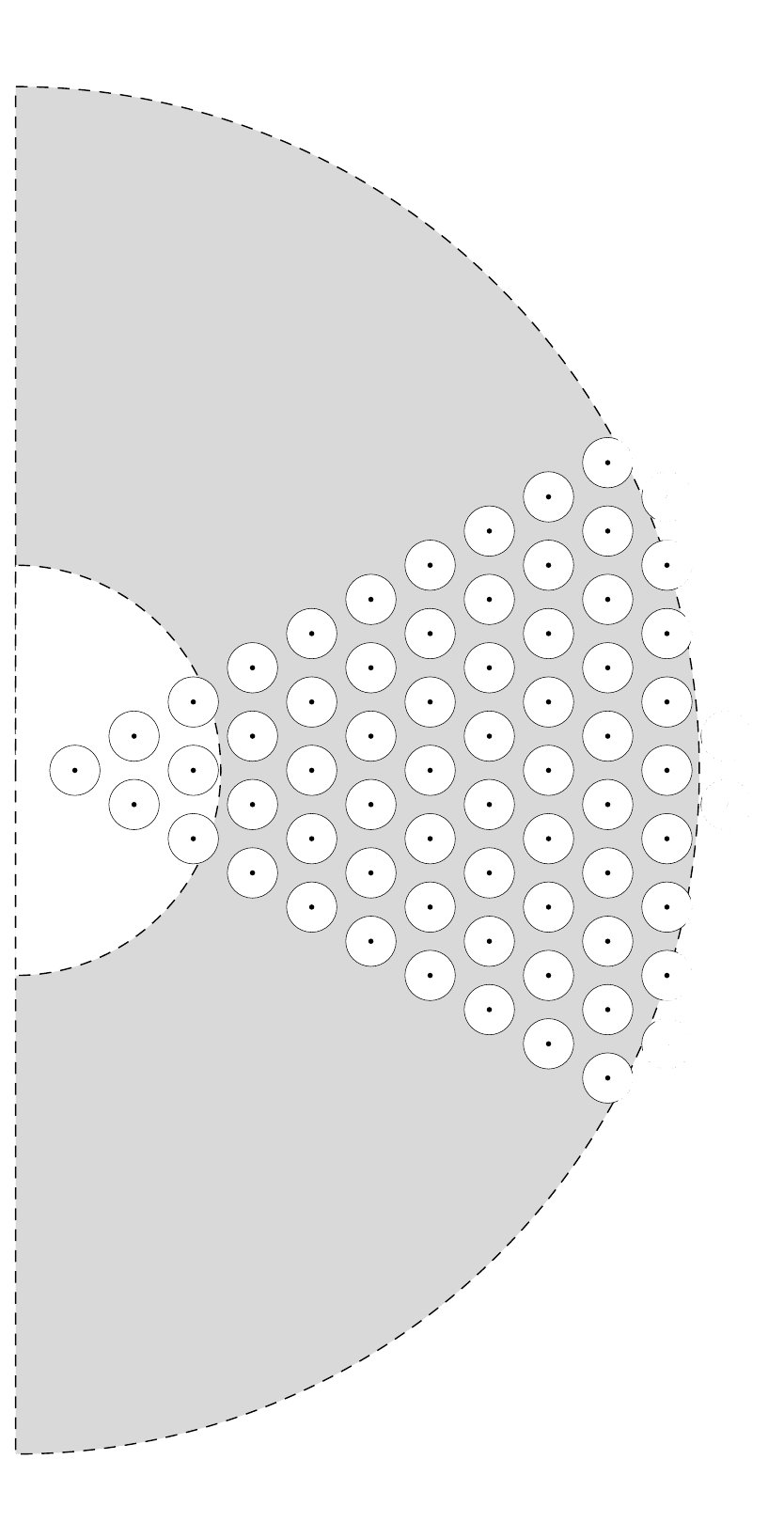}
\label{fig:figure2.1}
\end{minipage}
\hspace{0.5cm}
\end{figure}

\noindent
The disks surrounding the spectral values of the quadratic operator $q^w(x,D_x)$ correspond to the forbidden region 
$$\textrm{dist}(z,\mathrm{Spec}(q^w(x,D_x))) < h e^{-\frac{1}{C_1}(\log \log \frac{1}{h})^{1/n}}.$$
By coming back to the resolvent estimate (\ref{ning4}), we notice that the estimate
\begin{equation}\label{ning4123}
 h^{\frac{2k_0}{2k_0+1}}\|u\|_{L^2} \leq C_0 \|Pu-zu\|_{L^2}, 
\end{equation}
holds true at the boundary of the parabolic set $\Omega_{h}$, when 
$$\textrm{Re }z \leq c_1h^{\frac{2k_0}{2k_0+1}}, \quad \Big||\textrm{Im }z|-\frac{c_0}{2}\Big| \leq c_1,$$
with $0<c_1 \ll 1$. By using semigroup techniques, this resolvent estimate was improved by Sj\"ostrand~\cite{sjosemi} to 
\begin{equation}\label{ning41234b}
|\textrm{Re }z|\|u\|_{L^2} \leq C_0\|Pu-zu\|_{L^2}, 
\end{equation}
when 
$$-c_1 \leq \textrm{Re }z \leq -h^{\frac{2k_0}{2k_0+1}}, \quad \Big||\textrm{Im }z|-\frac{c_0}{2}\Big| \leq c_1,$$
and
\begin{equation}\label{ning41234}
 h^{\frac{2k_0}{2k_0+1}}\|u\|_{L^2} \leq C_0\exp\Big(\frac{C_0}{h}(\textrm{Re }z)_+^{\frac{2k_0+1}{2k_0}}\Big) \|Pu-zu\|_{L^2}, 
\end{equation}
when
$$-h^{\frac{2k_0}{2k_0+1}} \leq \textrm{Re }z \leq c_1\Big(h \log \frac{1}{h}\Big)^{\frac{2k_0}{2k_0+1}} , \quad \Big||\textrm{Im }z|-\frac{c_0}{2}\Big| \leq c_1.$$
For $\textrm{Re }z \sim h^{\frac{2k_0}{2k_0+1}}$, we recover the estimate (\ref{ning4123}). Furthermore, this result shows that the spectral parameter may enter 
logarithmically deeper into the numerical range outside the parabolic region~$\Omega_h$, 
$$\textrm{Re }z \sim\Big(h \log \frac{1}{h}\Big)^{\frac{2k_0}{2k_0+1}},$$
while keeping a polynomial resolvent bound 
$$\|(P-z)^{-1}\|_{\mathscr{L}(L^2)}=\mathcal{O}(h^{-\frac{2k_0}{2k_0+1}-\rho_0}), \quad \rho_0>0.$$

In the present work, we aim at completing this picture by investigating further the pseudospectral properties of the operator $P$ inside the $h$-neighborhood of the set
$$\Sigma=p_1(0)+\mathrm{Spec}(q^w(x,D_x)).$$
More specifically, we study necessary and sufficient conditions for the following a priori estimates to hold:  
\begin{multline}\label{sl2}
\exists c_0>0, \exists 0<h_0 \leq 1, \forall u \in \mathscr{S}(\rr^n), \forall \ 0 <h \leq h_0, \\ \|Pu-hzu\|_{L^2} \geq c_0 h^{\frac{N_0}{2}+1}\|u\|_{L^2},
\end{multline}
where $N_0 \geq 1$ is a positive integer, when $z$ belongs to a neighborhood of~$\Sigma$. While the resolvent estimate (\ref{ning4}) and the geometry of 
the parabolic region (\ref{sl1}) were shown to be related to the subelliptic properties of the quadratic approximation $q^w(x,D_x)$, we show in this work that 
the resolvent estimates (\ref{sl2}) are actually linked to some properties of hypoellipticity with a big loss of derivatives.  
The proof of the main result of this article (Theorem~\ref{th1}) is indeed based on a Grushin-reduction method following closely and adapting to the semiclassical 
setting the approach developed by Parenti and the first author in the study of hypoellipticity with a big loss of derivatives for operators with multiple characteristics~\cite{PP1}.

We recall that the Grushin-reduction method has proven itself fundamental in many different problems, especially in spectral theory and in the study of hypoellipticity (and, more recently in the study of
solvability and semi-global solvability \cite{S, He, PP1, PP2, PP3}) for operators with multiple characteristics. The idea of Grushin is roughly the following. 
Let $\mathsf{H}$ be a Hilbert space and $A\colon\mathsf{H}\longrightarrow\mathsf{H}$ be a Fredholm operator (with a non-zero kernel) of index $0$. 
Setting $V_1=\textrm{Ker }A$ and $V_2=\textrm{Ker }A^*$, these vector subspaces must then have the same dimension $d<+\infty$. 
We consider some orthonormal bases on these two vector subspaces $V_1=\mathrm{Span}\{w_1,\ldots,w_d\}$ and $V_2=\mathrm{Span}\{v_1,\ldots,v_d\}$.
Given a vector subspace $V\subset\mathsf{H}$ spanned by the orthonormal basis $\{e_1,\ldots,e_d\}$, we define the maps
$$h^+_V\colon \mathsf{H}\ni u\longmapsto\left[\begin{array}{c}(u,e_1)_\mathsf{H}\\
\vdots\\(u,e_d)_\mathsf{H}\end{array}\right]\in\mathbb{C}^d,\,\,\,\,h^-_V\colon\mathbb{C}^d\ni v=\left[\begin{array}{c}\zeta_1\\
\vdots\\\zeta_d\end{array}\right]\longmapsto\sum_{j=1}^d \zeta_je_j\in V.$$
Then, the system 
$$\mathbf{A}=\left[\begin{array}{cc}A & h^-_{\textrm{Ker }A^*}\\ h^+_{\textrm{Ker } A} & 0\end{array}\right]\colon \mathsf{H}\times\mathbb{C}^d\longrightarrow \mathsf{H}\times\mathbb{C}^d,$$
is invertible, with an inverse of the form
$$\mathbf{E}=\left[\begin{array}{cc}\tilde{E} & h^-_{\textrm{Ker }A}\\ h^+_{\textrm{Ker }A^*} & 0\end{array}\right]
\colon \mathsf{H}\times\mathbb{C}^d\longrightarrow \mathsf{H}\times\mathbb{C}^d,$$
where 
\begin{align*}
\tilde{E} : H=V_2 \oplus V_2^{\perp} &\rightarrow V_1^{\perp}\\
u=u_1+u_2 &\mapsto (A|_{V_1^{\perp}})^{-1}u_2.
\end{align*}
Next, given an operator $P$ with multiple characteristics, we take for $A$ a polynomial-coefficient differential operator
called \textit{the localized operator}~\cite{B,S}. The latter is the Weyl-quantization, in the normal directions to the characteristic set, of  
the relevant piece in the Taylor expansion of the symbol at the characteristic points obtained by keeping track of the vanishing orders of its various parts. Then, the system 
$$\left[\begin{array}{cc}P & R_-\\ R_+ & 0\end{array}\right],$$ which is approximated by the
system $\mathbf{A}$, can be inverted in a suitable pseudodifferential calculus (in the sense of left and right parametrices) by a system  
$$\left[\begin{array}{cc} E & E_-\\ E_+ & E_\pm\end{array}\right],$$ 
which is approximated by the system $\mathbf{E}$.

As already mentioned, this method has proven itself to be successful in the paper by Sj\"ostrand \cite{S} in which, for
the first time, the hypoellipticity with a loss of one derivative (and solvability) for general pseudodifferential operators with multiple characteristics was studied. Afterwards,
in this line of problems, Helffer \cite{He} studied the hypoellipticity with a loss of $3/2$-derivatives for operators with multiple characteristics. 
Pushing the machinery of localized operators to all orders (to describe the ``transport terms'' in the parametrix), 
Parenti and the first author~\cite{PP1} studied the hypoellipticity
with a big loss of derivatives for operators with multiple symplectic characteristics. They showed in particular that the various examples of $C^\infty$ hypoelliptic operators
with multiple characteristics and loss of derivatives, such as the Stein example, the Christ flat-Kohn example and others, where manifestations
of the same phenomenon~\cite{PP2}. More recently, they could also obtain, by the approach developed in \cite{PP1}, the local and semi-global solvability of certain
operators with multiple symplectic characteristics~\cite{PP3}.

We close this section by giving the plan of the article. The next section provides the statement of the main result (Theorem~\ref{th1}). Section~\ref{cases} is dedicated to some study cases, 
whereas the proof of Theorem~\ref{th1} is given in Section~\ref{proof}. 
Finally, an Appendix gathering miscellaneous facts and notations about quadratic differential operators is provided in Section~\ref{appendix}.

\section{Statement of the main result}\label{statement}
We consider the semiclassical pseudodifferential operator $P$ given in
(\ref{sl4}) whose Weyl symbol $p(x,\xi;h)$ admits the semiclassical asymptotic expansion (\ref{sl5}) in the symbol class $S(1)$, and 
we assume that the principal symbol $p_0$ satisfies the assumptions (\ref{sl6}), (\ref{sl7}), (\ref{sl8}), (\ref{sl9}).

Let $q$ be the quadratic term in the Taylor expansion of the principal symbol at the doubly characteristic point $X=0$,
\begin{equation}\label{sl11}
p_0(X)=q(X)+\mathcal{O}(X^3), \quad X=(x,\xi) \in \rr^{2n},
\end{equation}
when $X\rightarrow 0$. The assumption $\textrm{Re }p_0 \geq 0$ implies that the complex-valued quadratic form $q$ has also a nonnegative real part,
$\textrm{Re }q\geq 0$.

In the present work, we do not consider the degenerate case, that is the case
when the quadratic form $q$ is only partially elliptic (i.e. it satisfies the ellipticity condition (\ref{kps223}) on its singular space).
Indeed, we assume the quadratic form $q$ to be \textit{elliptic} on the whole phase space 
\begin{equation}\label{sl12}
(x,\xi) \in \rr^{2n}, \quad q(x,\xi)=0\, \Longrightarrow\,  (x,\xi)=0.
\end{equation}
Under this assumption, the spectrum of the quadratic operator 
$$q^w(x,D_x)u(x)=\frac{1}{(2\pi)^n}\int_{\rr^{2n}}e^{i(x-y) \cdot \xi}q\Big(\frac{x+y}{2},\xi\Big)u(y)d\xi dy,$$
is only composed of eigenvalues with finite algebraic multiplicities~\cite{S} (Theorem~3.5),
\begin{equation}\label{sl60}
\mathrm{Spec}(q^w(x,D_x))=\Big\{ \sum_{\substack{\lambda \in \mathrm{Spec}(F), \\  -i \lambda \in \Sigma(q)} }
{\big{(}r_{\lambda}+2 k_{\lambda}
\big{)}(-i\lambda) ; k_{\lambda} \in \nn }
\Big\},
\end{equation}
where $\Sigma(q)=\overline{q(\rr^{2n})}$, and
where $r_{\lambda}$ is the dimension of the complex vector space spanned by the generalized eigenvectors associated with the eigenvalue 
$\lambda \in \cc$ of the Hamilton map of $q$, see the Appendix (Section~\ref{appendix}).

Let $K \subset \cc$ be a compact set and let $N_0 \geq 1$ be a positive integer. We consider a spectral parameter $z(h)$
with the following semiclassical expansion
\begin{equation}\label{sl13}
z(h)=\sum_{k=0}^{2N_0+2}z_k h^{k/2},
\end{equation}
with $z_k \in K$ for all $0 \leq k \leq 2N_0+2$, where the leading term is assumed to satisfy
\begin{equation}\label{sl14}
z_0 \in p_1(0)+\mathrm{Spec}(q^w(x,D_x)).
\end{equation}
We define the symbols
\begin{equation}\label{sl15}
a_k(X)=\tilde{a}_k(X)-z_k:=\sum_{\substack{ j+\frac{|\alpha|}{2}=1+\frac{k}{2}\\ 0 \leq j \leq 1+[\frac{N_0}{2}], \ |\alpha| \leq N_0+2}}\frac{p_j^{(\alpha)}(0)}{\alpha!}X^{\alpha}-z_k,
\end{equation}
for $0 \leq k \leq 2N_0+2$, where $[x]$ stands for the integer part of $x$. Notice that
\begin{equation}\label{sl16}
a_0(X)=q(X)+p_1(0)-z_0,
\end{equation}
with $q$ the quadratic form defined in (\ref{sl11}).
The two operators 
\begin{equation}\label{sl17}
Q=a_0^w(x,D_x)=q^w(x,D_x)+p_1(0)-z_0 : B \longrightarrow L^2(\rr^n), 
\end{equation}
\begin{equation}\label{sl18}
Q^*=\overline{a_0}^w(x,D_x)=\overline{q}^w(x,D_x)+\overline{p_1(0)}-\overline{z_0} : B \longrightarrow L^2(\rr^n),  
\end{equation}
are known to be Fredholm operators of index 0 (see~\cite{H} , Lemma~3.1, or~\cite{S}, Theorem~3.5),
where $B$ is the Hilbert space
$$B=\{u \in L^2(\rr^n) ; x^{\alpha}D_x^{\beta}u \in L^2(\rr^n), \ \alpha,\beta\in \nn^n, \ |\alpha+\beta|\leq 2 \},$$
equipped with the norm
$$\|u\|_B^2=\sum_{|\alpha+\beta| \leq 2}\|x^{\alpha}D_x^{\beta}u\|_{L^2}^2.$$
Setting
\begin{equation}\label{sl19}
V_1=\textrm{Ker }Q, \quad  V_2=\textrm{Ker }Q^*,
\end{equation}
we may decompose $L^2(\rr^n)$ as
\begin{equation}\label{sl20}
L^2(\rr^n)=V_1 \oplus V_1^{\perp}=V_2 \oplus V_2^{\perp},
\end{equation}
with $V_1^{\perp}=\textrm{Ran }Q^*$, $V_2^{\perp}=\textrm{Ran }Q$. Since
$$0=\textrm{ind }Q=\textrm{dim Ker }Q-\textrm{codim Ran }Q=\textrm{dim }V_1-\textrm{codim }V_2^{\perp},$$
the kernels of the two operators $Q$, $Q^*$ have the same dimension 
$$1 \leq d=\textrm{dim }V_1=\textrm{dim }V_2 <+\infty.$$
Let $\phi_1,...,\phi_d$, respectively $\psi_1,...,\psi_d$, be an orthonormal basis of $V_1$, respectively $V_2$, whence
\begin{equation}\label{sl21}
Q\phi_j=0, \quad Q^*\psi_k=0, \quad 1 \leq j,k \leq d.
\end{equation} 
Because of the ellipticity of the quadratic symbols $q$ and $\overline{q}$, the eigenfunctions $\phi_j$, $\psi_k$ belong to the Schwartz space $\mathscr{S}(\rr^n)$. 
We denote respectively by $\pi_1$ and $\pi_2$ the orthogonal projections onto the vector spaces $V_1^{\perp}$ and $V_2^{\perp}$,
\begin{equation}\label{sl22}
\pi_1u=u-\sum_{j=1}^d(u,\phi_j)_{L^2}\phi_j, \quad \pi_2u=u-\sum_{j=1}^d(u,\psi_j)_{L^2}\psi_j.
\end{equation}
The two unbounded operators
\begin{equation}\label{sl80}
Q|_{V_1^{\perp}} : V_1^{\perp} \longrightarrow V_2^{\perp}=\textrm{Ran }Q, \quad Q^*|_{V_2^{\perp}} : V_2^{\perp} \longrightarrow V_1^{\perp}=\textrm{Ran }Q^*,
\end{equation}
are isomorphisms when equipped with the domains
$$D(Q|_{V_1^{\perp}})=B \cap V_1^{\perp}, \quad D(Q^*|_{V_2^{\perp}})=B \cap V_2^{\perp}.$$
We define the operator
\begin{align}\label{sl23}
S :  L^2(\rr^n)&=V_2 \oplus V_2^{\perp} \longrightarrow L^2(\rr^n) \\
 u&=u_1+u_2  \mapsto (Q|_{V_1^{\perp}})^{-1}u_2. \notag
\end{align}
The main result contained in this article is hence given by the following theorem.

\bigskip

\begin{theorem}\label{th1}
Let $K \subset \cc$ be a compact subset and $N_0 \geq 1$ be a positive integer. 
Let $P$ be a semiclassical pseudodifferential operator (\ref{sl4}) satisfying the assumptions (\ref{sl5}), (\ref{sl6}), (\ref{sl7}), (\ref{sl8}), (\ref{sl9}), 
(\ref{sl12}). Let $z(h)$ be the spectral parameter (\ref{sl13}) whose leading part satisfies to the assumption (\ref{sl14}). Let $\Omega$ be a compact subset of $K^{2N_0+2}$
(the Cartesian product of $K$ with itself $2N_0+2$ times).
The a priori estimate
\begin{multline}\label{sl31}
\exists c_0>0,\exists 0<h_0 \leq 1, \forall u \in L^2(\rr^n), \forall 0 <h \leq h_0, \\
\forall (z_1,...,z_{2N_0+2}) \in \Omega, \quad  \|Pu-hz(h)u\|_{L^2} \geq c_0 h^{\frac{N_0}{2}+1}\|u\|_{L^2},
\end{multline}
holds if and only if the a priori estimate 
\begin{multline}\label{sl32}
\exists c_0>0,\exists 0<h_0 \leq 1, \forall u_- \in \cc^d, \forall 0 <h \leq h_0, \\
\forall (z_1,...,z_{2N_0+2}) \in \Omega, \quad  |E_{\pm}u_-| \geq c_0 h^{\frac{N_0}{2}+1}|u_-|,
\end{multline}
holds, where $E_{\pm}$ stands for the $d\times d$ matrix 
\begin{equation}\label{sl30}
E_{\pm}=\sum_{j=1}^{2N_0+2}A_jh^{1+\frac{j}{2}}, \quad A_j=\big(A_{k,l}^{(j)}\big)_{1 \leq k,l \leq d} \in M_d(\cc),
\end{equation}
and where the entries $A_{k,l}^{(j)}$ of each $A_j$ are given by
\begin{equation}\label{sl33}
A_{k,l}^{(j)}=\sum_{i=1}^j(-1)^{i}\sum_{\substack{1 \leq k_p \leq 2N_0+2 \\ k_1+...+k_i=j}} (a_{k_1}^wS a_{k_2}^wS\ ...\ a_{k_{i-1}}^wS a_{k_i}^w\phi_{l},\psi_k)_{L^2}.
\end{equation}
The operators $a_k^w(x,D_x)$ are the Weyl quantizations of the symbols defined in~(\ref{sl15}).
\end{theorem}

\bigskip

\noindent
\begin{remark}\label{rem} It will be shown in the proof of Theorem~\ref{th1} that the operator $S$ is a pseudodifferential operator. 
This accounts for the fact that $S : \mathscr{S}(\rr^n) \longrightarrow \mathscr{S}(\rr^n)$ and the definition of the entries $A_{k,l}^{(j)}$. 
\end{remark}

\section{Some study cases}\label{cases}

Before plunging into the proof of Theorem \ref{th1}, which will be given in Section~\ref{proof}, we wish
to discuss in this section some study cases. We begin by studying the case of semiclassical hypoelliptic estimates with a loss of $3/2$ derivatives.

\subsection{Semiclassical hypoelliptic estimates with a loss of $\mathbf{3/2}$ derivatives}
When $N_0=1$, Theorem~\ref{th1} shows that the semiclassical hypoelliptic estimate with a loss of $3/2$ derivatives
\begin{multline}\label{sl41}
\exists c_0>0,\exists 0<h_0 \leq 1, \forall u \in L^2(\rr^n), \forall 0 <h \leq h_0, \forall  z_1 \in K, \\ \|Pu-hz_0u-h^{3/2}z_1u\|_{L^2} \geq c_0 h^{\frac{3}{2}}\|u\|_{L^2},
\end{multline}
holds, if and only if the matrix 
$$A_1(z_1)=((a_{1}^w(x,D_x)\phi_{l},\psi_k)_{L^2})_{1 \leq k,l \leq d},$$ 
is invertible for all $z_1 \in K$, where $a_1^w(x,D_x)$ is the differential operator defined by the Weyl quantization of the symbol
$$a_1(X)=\sum_{|\alpha|=3}\frac{p_0^{(\alpha)}(0)}{\alpha!}X^{\alpha}+\sum_{|\alpha|=1}\frac{p_1^{(\alpha)}(0)}{\alpha!}X^{\alpha}-z_1.$$

Denoting by $d_0$ the rank of the matrix $A_0=((\phi_{l},\psi_k)_{L^2})_{1 \leq k,l \leq d},$ we notice that its determinant $\det A_1(z_1)$ is a
polynomial function in the variable $z_1$ of degree~$d_0$. We distinguish two cases:

\medskip

\begin{itemize}
\item[(i)] $d_0=0$;
\item[(ii)] $1\leq d_0\leq d$.
\end{itemize}

\medskip

When $d_0=0$, the matrix $A_0$ is zero and the invertibility of the matrix $A_1(z_1)=A_1(0)$ is independent of the parameter $z_1$. When $\det A_1(0) \neq 0$, the a priori estimate (\ref{sl41}) holds. 
This indicates that there is no eigenvalue for the operator $P$ in any $h^{3/2}$-neighborhood of the point $hz_0$, when $0<h \ll 1$. On the other hand, when $\det A_1(0)=0$, the a priori estimate (\ref{sl41}) 
is violated for every $z_1 \in \cc$ and the resolvent cannot be bounded in norm as $\mathcal{O}(h^{-3/2})$ in any $h^{3/2}$-neighborhood of the point $hz_0$.

\medskip
When $1\leq d_0 \leq d$, we consider an open neighborhood $\omega$ of the finite set
$$\Lambda=\{z \in \cc ; \det A_1(z)=0\}.$$  
We deduce from Theorem~\ref{th1} that 
\begin{multline*}
\exists c_0>0, \exists 0<h_0 \leq 1, \forall u \in \mathscr{S}(\rr^n), \forall \ 0 <h \leq h_0, \forall z_1 \in K \cap (\cc \setminus \omega), \\ \|Pu-hz_0u-h^{3/2}z_1u\|_{L^2} \geq c_0 h^{\frac{3}{2}}\|u\|_{L^2}.
\end{multline*}
In this case, the spectrum of the operator $P$ in the disk $D(hz_0,Ch^{3/2})$
is localized in any $h^{3/2}$-neighborhood $U$
of the set $hz_0+h^{3/2}\Lambda$, and the resolvent of $P$ is bounded in norm as $\mathcal{O}(h^{-3/2})$ on the set $D(hz_0,Ch^{3/2}) \cap (\cc \setminus U)$.

\subsection{Case when the eigenfunctions have some parity properties}
When the eigenfunctions  $\phi_1,...,\phi_d$, $\psi_1,...,\psi_d$ enjoy some parity properties, the conclusions of Theorem~\ref{th1} may be 
further specified as follows.

\bigskip

\begin{proposition}\label{prop3}
Under the hypotheses of Theorem~\ref{th1}, we make the additional assumptions:  
\begin{itemize}
\item[$(i)$] The functions $\phi_1,...,\phi_d$ are all even, or all odd
\item[$(ii)$] The functions $\psi_1,...,\psi_d$ are all even, or all odd
\item[$(iii)$] All the terms with odd indices in the semiclassical expansion of the spectral parameter (\ref{sl13}) are zero
$$z(h)=\sum_{k=0}^{N_0+1}z_{2k} h^{k}$$
\end{itemize}
Then, the conclusions of Theorem~\ref{th1} holds with 
$$E_{\pm}=\sum_{j=1}^{2N_0+2}A_jh^{1+\frac{j}{2}},$$ 
where
$$A_{2j+1}=0,\quad\forall j\,\,\mathit{with}\,\,\ 1 \leq 2j+1 \leq 2N_0+2,$$
when the functions $\phi_1,...,\phi_d$ and $\psi_1,...,\psi_d$ have the same parity,  
or else
$$A_{2j}=0,\quad\forall j\,\,\mathit{with}\,\,\ 1 \leq 2j \leq 2N_0+2,$$
when the functions $\phi_1,...,\phi_d$ and $\psi_1,...,\psi_d$ have opposite parity. 
\end{proposition}

\bigskip

\begin{proof}
To begin with, we claim that 
\textit{when the functions $\psi_1,...,\psi_d \in V_2$ are all even, or all odd, then $Su$ is even (respectively odd) whenever $u \in \mathscr{S}(\rr^n)$ is an even (respectively odd) function.}\\
To see this, observe that when all the functions $\psi_1,...,\psi_d$ are even, the function
$$\pi_2 u=u-\sum_{j=1}^d(u,\psi_j)_{L^2}\psi_j,$$ 
is even (respectively odd) whenever $u \in \mathscr{S}(\rr^n)$ is even (respectively odd), because $(u,\psi_j)_{L^2}=0$ when $u$ is odd. On the other hand, when all the functions $\psi_1,...,\psi_d$ 
are odd, $\pi_2 u$ is also even (respectively odd) whenever $u \in \mathscr{S}(\rr^n)$ is even (respectively odd) because $(u,\psi_j)_{L^2}=0$ when $u$ is even. Then, we notice that $Qu$ is even 
(respectively odd) whenever $u \in \mathscr{S}(\rr^n)$ is even (respectively odd). Indeed, recalling that 
$$Q=q^w(x,D_x)+p_1(0)-z_0,$$ 
the parity property holds true for $Qu$ since it trivially holds true in the case of the operators
$$(x^{\alpha}\xi^{\beta})^w=\frac{1}{2}(x^{\alpha}D_x^{\beta}+D_x^{\beta}x^{\alpha}), \quad |\alpha+\beta|=2.$$
For $u \in \mathscr{S}(\rr^n)$, we write $Su=v_1+v_2$ with $v_1,v_2 \in \mathscr{S}(\rr^n)$ with $v_1$ even and $v_2$ odd (see Remark~\ref{rem}).
We assume that $u \in \mathscr{S}(\rr^n)$ is even (respectively odd). It follows from (\ref{sl80}) and (\ref{sl23}) that 
$$\textrm{Ran }S \subset V_1^{\perp}, \quad \pi_2=QS.$$
Since $\pi_2u=QSu=Qv_1+Qv_2$ is even (respectively odd), then $Qv_2=0$ (respectively $Qv_1=0$), that is, $v_2 \in V_1$ 
(respectively $v_1 \in V_1$). On the other hand, we have
$$0=(Su,v_2)_{L^2}=(v_1+v_2,v_2)_{L^2}=(v_1,v_2)_{L^2}+\|v_2\|_{L^2}^2=\|v_2\|_{L^2}^2,$$
respectively 
$$0=(Su,v_1)_{L^2}=(v_1+v_2,v_1)_{L^2}=\|v_1\|_{L^2}^2+(v_2,v_1)_{L^2}=\|v_1\|_{L^2}^2,$$
because $\textrm{Ran }S \subset V_1^{\perp}$, and $(v_1,v_2)_{L^2}=0$ when $v_1,v_2$ have opposite parity. It follows that $Su=v_1$ is even (respectively $Su=v_2$ is odd). This ends the proof of the claim. 

Next,  we notice that 
\begin{multline*}
\big(a^w(x,D_x)u\big)(-x)=\frac{1}{(2\pi)^n}\int_{\rr^{2n}}e^{i(-x-y)\cdot \xi}a\Big(\frac{-x+y}{2},\xi\Big)u(y)dyd\xi\\ 
=\frac{1}{(2\pi)^n}\int_{\rr^{2n}}e^{i(x-y)\cdot \xi}a\Big(-\frac{x+y}{2},-\xi\Big)u(-y)dyd\xi.
\end{multline*}
It follows that the function $a^wu$ is even (respectively odd) whenever $u \in \mathscr{S}(\rr^n)$ is even (respectively odd) when the symbol $a$ is even, whereas $a^wu$ is odd (respectively even) 
whenever $u \in \mathscr{S}(\rr^n)$ is even (respectively odd) when $a$ is odd. Therefore, when all the terms with odd indices in the semiclassical expansion of the spectral parameter (\ref{sl13}) are zero, that is
$$z(h)=\sum_{k=0}^{N_0+1}z_{2k} h^{k},$$
we have from (\ref{sl15}) that 
$$a_{2k}(X)=\sum_{\substack{ j+\frac{|\alpha|}{2}=1+k\\ 0 \leq j \leq 1+[\frac{N_0}{2}], \ |\alpha| \leq N_0+2}}\frac{p_j^{(\alpha)}(0)}{\alpha!}X^{\alpha}-z_{2k},$$
is an even function and that 
$$a_{2k+1}(X)=\sum_{\substack{ j+\frac{|\alpha|}{2}=1+k+\frac{1}{2}\\ 0 \leq j \leq 1+[\frac{N_0}{2}], \ |\alpha| \leq N_0+2}}\frac{p_j^{(\alpha)}(0)}{\alpha!}X^{\alpha},$$
is an odd function. Under assumption $(ii)$, we deduce from Remark~\ref{rem} and the previous claim that the function 
$$a_{k_1}^wS a_{k_2}^wS\ ...\ a_{k_{i-1}}^wS a_{k_i}^wu \in \mathscr{S}(\rr^n),\,\,\,\,\mathrm{with}\,\,\,\,1 \leq k_1+...+k_i=2j \leq 2N_0+2,$$
is even (respectively odd) whenever $u \in \mathscr{S}(\rr^n)$ is even (respectively odd). On the other hand, the function 
$$a_{k_1}^wS a_{k_2}^wS\ ...\ a_{k_{i-1}}^wS a_{k_i}^wu \in \mathscr{S}(\rr^n),\,\,\,\,\mathrm{with} \,\,\,\,1 \leq k_1+...+k_i=2j+1 \leq 2N_0+2,$$
is odd (respectively even) whenever $u \in \mathscr{S}(\rr^n)$ is even (respectively odd). It follows that 
$$(a_{k_1}^wS a_{k_2}^wS\ ...\ a_{k_{i-1}}^wS a_{k_i}^w\phi_{l},\psi_k)_{L^2}=0,$$
with $1 \leq k_1+...+k_i=2j+1 \leq 2N_0+2$ (respectively $1 \leq k_1+...+k_i=2j \leq 2N_0+2$), when the functions $\phi_1,...,\phi_d$ and $\psi_1,...,\psi_d$ 
have the same parity (respectively opposite parity). This ends the proof of Proposition~\ref{prop3}.
\end{proof}

\subsection{Case $\mathbf{d=\textrm{dim }V_1=\textrm{dim }V_2=1}$} We now consider the case when the kernels $V_1$ and $V_2$ are one-dimensional, that is
$$V_1=\textrm{Ker }Q=\textrm{Span }\phi_1, \quad V_2=\textrm{Ker }Q^*=\textrm{Span }\psi_1,$$
spanned by eigenfunctions satisfying $(\phi_1,\psi_1)_{L^2} \neq 0$.
In this case, the matrix (\ref{sl30}) can be written as
$$E_{\pm}=\sum_{j=1}^{2N_0+2}h^{1+\frac{j}{2}}\sum_{i=1}^j(-1)^{i}\sum_{\substack{1 \leq k_p \leq 2N_0+2 \\ k_1+...+k_i=j}} (a_{k_1}^wS a_{k_2}^wS\ ...\ a_{k_{i-1}}^wS a_{k_i}^w\phi_{1},\psi_1)_{L^2}.$$

We define successively for every $1 \leq j \leq 2N_0+2$,
\begin{equation}\label{eq17.1}
\tilde{z}_1=\frac{1}{(\phi_1,\psi_1)_{L^2}}(\tilde{a}_{1}^w(x,D_x)\phi_1,\psi_1)_{L^2},
\end{equation}
\begin{multline}\label{eq15.1}
\tilde{z}_j=\frac{1}{(\phi_1,\psi_1)_{L^2}}\Big[(\tilde{a}_j^w\phi_1,\psi_1)_{L^2} \\ +\sum_{i=2}^j(-1)^{i+1}\hspace{-.7cm}
\sum_{\substack{1 \leq k_p \leq 2N_0+2 \\ 
k_1+...+k_i=j}}\hspace{-.4cm} \big((\tilde{a}_{k_1}^w-\tilde{z}_{k_1})S (\tilde{a}_{k_2}^w-\tilde{z}_{k_2})S\ ...\ (\tilde{a}_{k_{i-1}}^w-\tilde{z}_{k_{i-1}})S (\tilde{a}_{k_i}^w-\tilde{z}_{k_i})\phi_{1},\psi_1)_{L^2}\Big].
\end{multline}
The following result follows from Theorem~\ref{th1}.

\bigskip

\begin{corollary}\label{cor1}
Under the hypotheses of Theorem~\ref{th1}, we assume further that the two kernels $V_1=\textrm{Span }\phi_1$, $V_2=\textrm{Span }\psi_1$ are one-dimensional, spanned by eigenfunctions satisfying $(\phi_1,\psi_1)_{L^2} \neq 0$.
Let $N_0 \geq 1$ be a positive integer and let $K \subset \cc \setminus \{\tilde{z}_{N_0}\}$ 
be a compact subset, where the complex numbers  $\tilde{z}_j$, $1\leq j\leq 2N_0+2$, are defined in (\ref{eq17.1}), (\ref{eq15.1}). 
Then, there exist $c_0>0$, $0<h_0 \leq 1$ such that for all $u \in L^2(\rr^n)$, $0 <h \leq h_0$, and $z \in K$ one has
$$\Big\|Pu-hz_0u-\sum_{j=1}^{N_0-1}h^{1+\frac{j}{2}}\tilde{z}_ju-h^{1+\frac{N_0}{2}}zu\Big\|_{L^2} \geq c_0 h^{1+\frac{N_0}{2}}\|u\|_{L^2}.$$
\end{corollary}

\bigskip

Recalling (\ref{sl60}) and (\ref{sl14}), we shall now consider the specific case when 
\begin{equation}\label{eq16.1}
z_0-p_1(0)= \sum_{\substack{\lambda \in \mathrm{Spec}(F), \\  -i \lambda \in \Sigma(q) }}-i \lambda r_{\lambda},
\end{equation}
is the first eigenvalue in the bottom of the spectrum of the elliptic quadratic operator $q^w(x,D_x)$.
In this case, we recall from \cite{S} (see also~\cite{kps11}, Theorem~2.1) that the eigenvalue $z_0-p_1(0)$ has algebraic multiplicity 1 and that the eigenspace 
$$V_1=\textrm{Ker }Q=\cc \phi_1,$$
is spanned by a ground state of exponential type 
$$\phi_1(x)=e^{-a(x)} \in \mathscr{S}(\rr^n),$$ 
where $a$ is a complex-valued quadratic form whose real part is positive definite. The quadratic form $a$ is defined as
$$a(x)=-\frac{1}{2}i\langle x,B^+x \rangle, \quad x \in \rr^n, \quad \textrm{Im }B^+>0,$$
where $B^+$ is the symmetric matrix with positive definite imaginary part $\mathrm{Im}\,B^+$
defining the positive Lagrangian plane (see \cite{S}, Proposition~3.3)
$$V^+=\bigoplus_{\substack{\lambda \in \mathrm{Spec}(F) \\  -i \lambda \in \Sigma(q)}}V_{\lambda}=\{(x,B^+x) ; x \in \cc^n\},$$
where $V_{\lambda}$ is the space of the 
generalized eigenvectors associated with the eigenvalue $\lambda$ of the Hamilton map of~$q$. 
On the other hand, we notice that the adjoint operator is actually given by the quadratic operator
$$q^w(x,D_x)^*=\overline{q}^w(x,D_x),$$
whose symbol is the complex conjugate symbol of $q$. This quadratic symbol is also elliptic. It follows that $\overline{z}_0-\overline{p_1(0)}$ 
is the first eigenvalue in the bottom of the spectrum for the quadratic operator $q^w(x,D_x)^*$. This eigenvalue has therefore algebraic multiplicity 1 and its eigenspace is 
also spanned by a ground state of exponential type 
$$\psi_1(x)=e^{-\tilde{a}(x)} \in \mathscr{S}(\rr^n),$$ 
where $\tilde{a}$ is a complex-valued quadratic form whose real part is positive definite. 
Under these assumptions, we are therefore in the situation where $d=\textrm{dim }V_1=\textrm{dim }V_2=1$, with even eigenfunctions $\phi_1, \psi_1$,
$$V_1=\textrm{Span } \phi_1, \quad V_2=\textrm{Span } \psi_1.$$
These two eigenspaces are equal $V_1=V_2$, that is $\phi_1=\psi_1$, if and only if $V^+=\overline{V^-}$ (see~\cite{viola3}, Theorem~1.7), where
$$V^-=\bigoplus_{\substack{\lambda \in \mathrm{Spec}(F) \\  i \lambda \in \Sigma(q)}}V_{\lambda}.$$
We consider the case when 
\begin{equation}\label{sl61}
(\phi_1,\psi_1)_{L^2} \neq 0, 
\end{equation}
and $N_0=2\tilde{N}_0$, with $\tilde{N}_0 \geq 1$. 
Notice that (\ref{sl61}) means that $\pi_{V_2}\bigl|_{V_1}\colon V_1\longrightarrow V_2$ is invertible, where $\pi_{V_2}\colon L^2(\mathbb{R}^n)\longrightarrow L^2(\mathbb{R}^n)$
is the orthogonal projection onto $V_2$ (of course, the same holds for $\pi_{V_1}\bigl|_{V_2}\colon V_2\longrightarrow V_1$).

Under these assumptions, we define successively for every $1 \leq j \leq 2\tilde{N}_0+1$,
\begin{equation}\label{eq16.3}
\tilde{z}_{2}=\frac{1}{(\phi_1,\psi_1)_{L^2}}[(\tilde{a}_{2}^w\phi_1,\psi_1)_{L^2}-(\tilde{a}_{1}^wS \tilde{a}_{1}^w\phi_{1},\psi_1)_{L^2}],
\end{equation}
\begin{multline}\label{eq16.4}
\tilde{z}_{2j}=\frac{1}{(\phi_1,\psi_1)_{L^2}}\Big[(\tilde{a}_{2j}^w\phi_1,\psi_1)_{L^2} \\ +\sum_{i=2}^{2j}(-1)^{i+1}\sum_{\substack{1 \leq k_p \leq 2N_0+2 \\ 
k_1+...+k_i=2j}} (a_{k_1}^wSa_{k_2}^wS\ ...\ a_{k_{i-1}}^wS a_{k_i}^w\phi_{1},\psi_1)_{L^2}\Big].
\end{multline}
We therefore deduce from Proposition~\ref{prop3} and Theorem~\ref{th1} the following result.

\bigskip

\begin{corollary}\label{cor2}
Under the hypotheses of Theorem~\ref{th1}, we make the additional assumptions (\ref{eq16.1}), (\ref{sl61}).
Let $\tilde{N}_0 \geq 1$ be a positive integer and let $K \subset \cc \setminus \{\tilde{z}_{2\tilde{N}_0}\}$ be a 
compact subset, where the complex numbers  $\tilde{z}_{2j}$, $1\leq j\leq 2\tilde{N}_0+1$, 
are defined in (\ref{eq16.3}), (\ref{eq16.4}). Then, there exist $c_0>0$, $0<h_0 \leq 1$ such that for all $u \in\mathscr{S}(\rr^n)$, $0 <h \leq h_0$, and $z \in K$,
one has
$$\Big\|Pu-hz_0u-\sum_{j=1}^{\tilde{N}_0-1}h^{1+j}\tilde{z}_{2j}u-h^{\tilde{N}_0+1}zu\Big\|_{L^2} \geq c_0 h^{\tilde{N}_0+1}\|u\|_{L^2}.$$
\end{corollary}

\bigskip

On the other hand, we may calculate explicitly all the terms in the semiclassical expansion (\ref{eq1.9}) of an eigenvalue of $P$ with leading term $hz_0$, when the assumptions of Corollary~\ref{cor2} are satisfied. 
Indeed, the semiclassical expansion is given in this case by 
$$z_{k} \sim h (\lambda_{k} + p_1(0) + h \lambda_{k,1} + h^{2} \lambda_{k,2} +\ldots),$$
since $z_0-p_1(0)$ is an eigenvalue with algebraic multiplicity~1 of the quadratic operator $q^w(x,D_x)$. Then, we directly deduce from Proposition~\ref{prop3} and Theorem~\ref{th1} 
that the coefficients $\lambda_{k,j}$ must correspond to the terms  $\tilde{z}_{2j}$ given in Corollary~\ref{cor2}, that is
$$\lambda_{k,j}=\tilde{z}_{2j},\quad 1 \leq j \leq \tilde{N}_0.$$

\section{Proof of Theorem~\ref{th1}}\label{proof}

Let $K \subset \cc$ be a compact subset and $N_0 \geq 1$ be a positive integer. Let
$$P=p^w(x,hD_x;h)=\frac{1}{(2\pi)^n}\int_{{\rr}^{2n}}e^{i(x-y)\cdot\xi}p\Big(\frac{x+y}{2},h\xi;h\Big)u(y)dyd\xi,$$
be a semiclassical pseudodifferential operator satisfying to the assumptions (\ref{sl5}), (\ref{sl6}), (\ref{sl7}), (\ref{sl8}), (\ref{sl9}), (\ref{sl12}), and let
$$z(h)=\sum_{k=0}^{2N_0+2}z_k h^{k/2}, \quad z_k \in K,$$
be the spectral parameter whose leading part satisfies assumption~(\ref{sl14}).
After conjugating by the unitary operator 
$$T_h : L^2(\rr^n) \ni u(x) \longmapsto h^{n/4}u(h^{1/2}x) \in L^2(\rr^n),$$
it is sufficient to prove Theorem~\ref{th1} for the operator 
$$P=T_hp^w(x,hD_x;h)T_h^{-1}=p^w(h^{1/2}x,h^{1/2}D_x;h).$$
In the following, the standard notation $c=a \#^w b$ denotes the Weyl symbol of the operator obtained by composition
$$c^w(x,D_x)=a^w(x,D_x)b^w(x,D_x)=(\textrm{Op}^w a)(\textrm{Op}^w b),$$
with the standard normalization of the Weyl quantization
$$a^w(x,D_x)u(x)=(\textrm{Op}^w a)u(x)=\frac{1}{(2\pi)^n}\int_{\rr^{2n}}e^{i(x-y)\cdot \xi}a\Big(\frac{x+y}{2},\xi\Big)u(y)dyd\xi.$$
We refer the reader to the notations introduced in Section~\ref{statement} and begin by noticing that the orthogonal projections $\pi_1$ and $\pi_2$ onto the vector spaces $V_1^{\perp}$ and $V_2^{\perp}$ satisfy
$$\pi_1, \pi_2 : \mathscr{S}(\rr^n) \longrightarrow \mathscr{S}(\rr^n),$$
since the eigenfunctions $\phi_j$, $\psi_k$ belong to the Schwartz space $\mathscr{S}(\rr^n)$. 
Next, we observe that the $L^2$-adjoint of the bounded operator $S : L^2(\rr^n) \rightarrow L^2(\rr^n)$ defined in (\ref{sl23}) is given by  
\begin{align}\label{sl23bis}
S^* :  L^2(\rr^n)&=V_1 \oplus V_1^{\perp} \longrightarrow L^2(\rr^n) \\
 u&=u_1+u_2  \mapsto (Q^*|_{V_2^{\perp}})^{-1}u_2. \notag
\end{align}
By definition, these two operators are continuous on $L^2(\rr^n)$ and satisfy the identities
\begin{equation}\label{sl24}
SQ=Q^*S^*=1-(1-\pi_1), \quad QS=S^*Q^*=1-(1-\pi_2).
\end{equation}
By using standard notations for metrics on the phase space~\cite{hormander,lerner}, for $m \in\mathbb{R}$ we write 
\begin{multline}\label{sl71}
\mathbf{S}^m=S(\langle X\rangle^m,g)=\{a \in C^{\infty}(\rr^{2n},\cc); \forall \alpha \in \nn^{2n}, \exists C_{\alpha}>0,\\  \forall X \in \rr^{2n},
\quad  |\partial_X^{\alpha} a(X)| \leq C_{\alpha}\langle X \rangle^{m-|\alpha|}\}.
\end{multline}
for the class of ($h$-independent) \textit{global} pseudodifferential operators (after Shubin, see \cite{shubin}), where   
$g$ is the admissible, geodesically temperate metric (see, e.g., \cite{lerner}, Lemma~2.6.23) given by
$$g=\frac{|dX|^2}{\langle X\rangle^2}, \quad X=(x,\xi) \in \rr^{2n},$$
and put $\mathrm{Op}^w(\mathbf{S}^m)$ for the set of corresponding pseudodifferential operators associated with symbols in $\mathbf{S}^m.$
 
As global pseudodifferential operators in $\rr^n$, the operators $1-\pi_1$ and $1-\pi_2$ are smoothing since their symbols in the standard quantization, 
respectively given by
$$\sum_{j=1}^de^{-i x \cdot \xi}\phi_j(x)\overline{\widehat{\phi}_j(\xi)} \in \mathbf{S}^{-\infty}, 
\quad \sum_{j=1}^de^{-i x \cdot \xi}\psi_j(x)\overline{\widehat{\psi}_j(\xi)}\in \mathbf{S}^{-\infty},$$
belong to the Schwartz space $\mathscr{S}(\rr^{2n})$. This also implies that
$$1-\pi_1 \in \textrm{Op}^w(\textbf{S}^{-\infty}), \quad 1-\pi_2 \in \textrm{Op}^w(\textbf{S}^{-\infty}).$$
Setting
$$T_+u=\sum_{j=1}^d(u,\psi_j)_{L^2}\phi_j, \quad T_+u=\sum_{j=1}^d(u,\phi_j)_{L^2}\psi_j,$$
the very same arguments show that $T_{\pm} \in \textrm{Op}^w(\textbf{S}^{-\infty})$. Next, we deduce from (\ref{sl19}), (\ref{sl21}) and (\ref{sl80}) that the mapping
\begin{align*}
\Phi : B &\rightarrow L^2(\rr^n)=V_2 \oplus V_2^{\perp}\\
u &\mapsto \sum_{j=1}^d(u,\phi_j)_{L^2}\psi_j+Qu,
\end{align*}
is invertible with inverse given by
\begin{align*}
\Phi^{-1} :  L^2(\rr^n) \rightarrow&\ B \\
u \mapsto& \ \sum_{j=1}^d(u,\psi_j)_{L^2}\phi_j+Su.
\end{align*}
We may write $\Phi=a^w(x,D_x)$, with $a \in \textbf{S}^2$.
By referring to~\cite{lerner} (Section~2.6) for the definition of Sobolev spaces attached to a pseudodifferential calculus, we notice that $L^2(\rr^n)=H(1,g)$ and 
$B=H(\langle X\rangle^2,g)$.
Then we deduce from \cite{lerner}, Corollary~2.6.28, that 
$$\Phi^{-1} \in\textrm{Op}^w(\textbf{S}^{-2}).$$
Since $T_{+} \in \textrm{Op}^w(\textbf{S}^{-\infty})$, this implies that $S \in \textrm{Op}^w(\textbf{S}^{-2})$ and justifies Remark~\ref{rem}.

We now get to the proof of Theorem~\ref{th1}. We start by noticing that for $X=(x,\xi) \in \rr^{2n}$,
$$p(h^{1/2}X;h)-hz(h)=\sum_{k=0}^{1+[\frac{N_0}{2}]}p_k(h^{1/2}X)h^k-\sum_{k=0}^{2N_0+2}z_k h^{1+\frac{k}{2}} \textrm{ mod } S(h^{[\frac{N_0}{2}]+2}),$$
up to a symbol belonging to the class $S(h^{[\frac{N_0}{2}]+2})$, whose definition is given in (\ref{sl70}).
By Taylor expanding the symbols
$$p_k(X)= \sum_{|\alpha| \leq N_0+2}\frac{p_k^{(\alpha)}(0)}{\alpha!}X^{\alpha}+\sum_{|\alpha|=N_0+3}\frac{N_0+3}{\alpha!}X^{\alpha}\int_0^1(1-t)^{N_0+2}p_k^{(\alpha)}(tX)dt,$$
we obtain that  
\begin{multline}\label{eq10.0}
p(h^{1/2}X;h)-hz(h)=\sum_{\substack{k=0,...,1+[\frac{N_0}{2}] \\ |\alpha| \leq N_0+2}}\frac{p_k^{(\alpha)}(0)}{\alpha!}X^{\alpha}h^{k+\frac{|\alpha|}{2}}-\sum_{k=0}^{2N_0+2}z_k h^{1+\frac{k}{2}} \\
+\sum_{\substack{k=0,...,1+[\frac{N_0}{2}] \\ |\alpha|=N_0+3}}\frac{N_0+3}{\alpha!}X^{\alpha}h^{k+\frac{|\alpha|}{2}}\int_0^1(1-t)^{N_0+2}p_k^{(\alpha)}(th^{1/2}X)dt \quad \textrm{mod } S(h^{[\frac{N_0}{2}]+2}).
\end{multline}
By assumption (\ref{sl9}), the point $0 \in \rr^{2n}$ is doubly characteristic for the principal symbol, that is $p_0(0)=\nabla p_0(0)=0$. 
Setting 
$$R_{\alpha}(X;h)=\sum_{k=0,...,1+[\frac{N_0}{2}]}\frac{N_0+3}{\alpha!}h^{k}\int_0^1(1-t)^{N_0+2}p_k^{(\alpha)}(th^{1/2}X)dt,$$
we may thus write 
\begin{multline}\label{eq10.1}
p(h^{1/2}X;h)-hz(h)\\ =\sum_{k=0}^{2N_0+2}a_k(X)h^{1+\frac{k}{2}}
+h^{\frac{N_0}{2}+\frac{3}{2}}\sum_{|\alpha|=N_0+3}X^{\alpha}R_{\alpha}(X;h) \textrm{ mod } S(h^{[\frac{N_0}{2}]+2}),
\end{multline}
with the symbols $a_k$ defined in (\ref{sl15}). As the $p_k \in S(1)$, we readily see that
\begin{equation}\label{sl72}
R_{\alpha} \in S(1).
\end{equation} 
Following~\cite{PP1}, we now use a Grushin-reduction method. To this end, we define  
\begin{align}\label{sl73}
R_- :   \cc^d & \  \longrightarrow V_2\\ \notag
u_- & \   \longmapsto \sum_{j=1}^du_-(j)\psi_j, 
\end{align}
\begin{align}\label{sl74}
R_+ :   L^2(\rr^n) & \  \longrightarrow \cc^d\\ \notag
u & \   \longmapsto ((u,\phi_j)_{L^2})_{1 \leq j \leq d},
\end{align}
where the $\phi_j, \psi_k$ are the eigenfunctions defined in (\ref{sl21}). 
Setting 
\begin{equation}\label{sl75}
\phi_{0,k}^+=\phi_k, \quad \psi_{0,k}^-=\psi_k, \quad  k=1,...,d,
\end{equation}
we shall construct by induction functions $\phi_{j,k}^+$, $\psi_{j,k}^- \in \mathscr{S}(\rr^n)$, 
for $1 \leq k \leq d$, $1 \leq j \leq 2N_0+2$, and $d\times d$ complex matrices $A_j \in M_d(\cc)$ for $j=1,...,2N_0+2$, 
which are all independent of the semiclassical parameter and satisfy the equations
\begin{equation}\label{eq4}
R_+E_+=\textrm{Id}+\mathcal{O}_{\mathscr{L}(\cc^d)}(h^{1/2}),
\end{equation}
\begin{equation}\label{eq3}
\sum_{k=0}^{2N_0+2}h^{1+\frac{k}{2}} a_k^w(x,D_x)E_++R_{-}E_{\pm}=\mathcal{O}_{\mathscr{L}(\cc^d,L^2)}(h^{N_0+\frac{5}{2}}),
\end{equation}
\begin{equation}\label{eq2}
E_-\left(\sum_{k=0}^{2N_0+2}a_k^w(x,D_x)h^{1+\frac{k}{2}}\right)+E_{\pm}R_+=\mathcal{O}_{\mathscr{L}(L^2,\cc^d)}(h^{N_0+\frac{5}{2}}),
\end{equation}
\begin{equation}\label{eq1}
S a_0^w(x,D_x) +E_+R_+=\textrm{Id}+\mathcal{O}_{\mathscr{L}(L^2)}(h^{1/2}),
\end{equation}
where
\begin{align}\label{eq10.2}
E_+ :   \cc^d & \  \longrightarrow L^2(\rr^n)\\
u_- & \   \longmapsto \sum_{k=1}^du_-(k)\left(\sum_{j=0}^{2N_0+2}\phi_{j,k}^+h^{\frac{j}{2}}\right), \notag
\end{align} 
\begin{align}\label{sl76}
E_- :   & \ L^2(\rr^n)   \longrightarrow \cc^d\\ \notag
& \ u  \longmapsto \Big(\Big(u,\sum_{j=0}^{2N_0+2}\psi_{j,k}^-h^{\frac{j}{2}}\Big)_{L^2}\Big)_{1 \leq k \leq d},
\end{align}
\begin{equation}\label{a1}
E_{\pm}=\sum_{j=1}^{2N_0+2}A_jh^{1+\frac{j}{2}}.
\end{equation}
The notation $\mathcal{O}_{\mathscr{L}(E,F)}(h^{N})$ stands for a remainder which is a bounded operator $T: E \rightarrow F$ with a norm satisfying $\|T\|_{\mathscr{L}(E,F)} \lesssim h^{N}$.

We next notice that equation (\ref{eq4}) is directly satisfied since the functions $(\phi_{0,k}^+)_{1 \leq k \leq d}=(\phi_k)_{1 \leq k \leq d}$ are chosen orthonormal
$$\Big(\sum_{k=1}^du_-(k)\Big(\sum_{j=0}^{2N_0+2}\phi_{j,k}^+h^{\frac{j}{2}}\Big),\phi_l\Big)_{L^2}=u_-(l)+\sum_{\substack{k=1,...,d \\ j=1,...,2N_0+2}}h^{\frac{j}{2}}u_-(k)(\phi_{j,k}^+,\phi_l)_{L^2}.$$
Equation (\ref{eq3}) can be written as
\begin{multline}\label{eq5}
\sum_{\substack{1 \leq l \leq d \\ 0 \leq j,k \leq 2N_0+2}}h^{1+\frac{k+j}{2}} u_-(l) a_k^w(x,D_x)\phi_{j,l}^+\\ +\sum_{\substack{1 \leq l \leq d \\ 1 \leq j \leq 2N_0+2}}h^{1+\frac{j}{2}}(A_ju_-)(l)\psi_l=\mathcal{O}(h^{N_0+\frac{5}{2}}|u_-|).
\end{multline}
We deduce from (\ref{sl17}), (\ref{sl21}) and (\ref{sl75}) that the term factoring $h$ in the left-hand side of (\ref{eq5}) is zero, that is,
$$\sum_{1 \leq l \leq d} u_-(l) a_0^w(x,D_x)\phi_{0,l}^+=\sum_{1 \leq l \leq d} u_-(l) Q\phi_{l}=0.$$
Next, we observe that the term factoring $h^{1+\frac{j}{2}}$, with $1 \leq j \leq 2N_0+2$, in the left-hand side of equation (\ref{eq5}) is zero if and only if we have
\begin{equation}\label{sl77}
\sum_{\substack{1 \leq l \leq d \\ 0 \leq k_1,k_2 \leq 2N_0+2 \\ k_1+k_2=j}} u_-(l) a_{k_1}^w(x,D_x)\phi_{k_2,l}^++\sum_{1 \leq l \leq d}(A_ju_-)(l)\psi_l=0.
\end{equation}
By assuming that the Schwartz functions $\phi_{k,l}^+$ have already been determined for all $0 \leq k \leq j-1$, $1 \leq l \leq d$, 
it will be sufficient for fulfilling equation (\ref{sl77}) to choose the functions $(\phi_{j,l}^{+})_{1 \leq l \leq d}$ and the matrix $A_j=(A_{k,l}^{(j)})_{1 \leq k,l \leq d}$ to satisfy the identities 
\begin{equation}\label{eq5.5}
Q\phi_{j,l}^+=-\sum_{1 \leq k \leq d}A_{k,l}^{(j)}\psi_k-\sum_{\substack{0 \leq k_1,k_2 \leq 2N_0+2 \\ k_1+k_2=j, \ k_1 \geq 1}} a_{k_1}^w(x,D_x)\phi_{k_2,l}^+,
\end{equation}
for every $1 \leq l \leq d$.
Taking
\begin{equation}\label{eq6}
A_{k,l}^{(j)}=-\sum_{\substack{0 \leq k_1,k_2 \leq 2N_0+2 \\ k_1+k_2=j, \ k_1 \geq 1}} (a_{k_1}^w(x,D_x)\phi_{k_2,l}^+,\psi_k)_{L^2},
\end{equation}
yields that the right-hand side of (\ref{eq5.5}) is fully determined and belongs to $V_2^{\perp}$,
$$Q\phi_{j,l}^+=-\pi_2\Big(\sum_{\substack{0 \leq k_1,k_2 \leq 2N_0+2 \\ k_1+k_2=j, \ k_1 \geq 1}} a_{k_1}^w(x,D_x)\phi_{k_2,l}^+\Big) \in V_2^{\perp}.$$
It follows from (\ref{sl80}) and (\ref{sl23}) that we can choose the functions $\phi_{j,l}^+$ as
\begin{align}\label{eq6.5}
\phi_{j,l}^+= & \ -(Q|_{V_1^{\perp}})^{-1}\pi_2\Big(\sum_{\substack{0 \leq k_1,k_2 \leq 2N_0+2 \\ k_1+k_2=j, \ k_1 \geq 1}} a_{k_1}^w(x,D_x)\phi_{k_2,l}^+\Big)\\
= & \ -S\Big(\sum_{\substack{0 \leq k_1,k_2 \leq 2N_0+2 \\ k_1+k_2=j, \ k_1 \geq 1}} a_{k_1}^w(x,D_x)\phi_{k_2,l}^+\Big). \notag
\end{align}
Since $S : \mathscr{S}(\rr^n) \longrightarrow \mathscr{S}(\rr^n)$ (see Remark~\ref{rem}), the functions $(\phi_{j,l}^+)_{1 \leq l \leq d}$ belong to $\mathscr{S}(\rr^n)$. 
By iterating this process, we obtain functions  $\phi_{j,l}^+$ for $1 \leq l \leq d$, $0 \leq j \leq 2N_0+2$ and matrices $A_j$ for $1 \leq j \leq 2N_0+2$, 
satisfying equation (\ref{eq3}). Next, equation (\ref{eq2}) can be written as 
\begin{multline*}
\sum_{0 \leq j,k \leq 2N_0+2}h^{1+\frac{j+k}{2}}(a_k^w(x,D_x)u,\psi_{j,l}^-)_{L^2}\\ +\sum_{\substack{1 \leq k \leq d \\ 1 \leq j \leq 2N_0+2}}h^{1+\frac{j}{2}}A_{l,k}^{(j)}(u,\phi_k)_{L^2}=\mathcal{O}(h^{N_0+\frac{5}{2}}\|u\|_{L^2}),
\end{multline*}
for $1 \leq l \leq d$. We need therefore to satisfy the equations
$$\Big(u,\hspace{-0.3cm}\sum_{0 \leq j,k \leq 2N_0+2}h^{1+\frac{j+k}{2}}\overline{a}_k^w(x,D_x)\psi_{j,l}^-+\hspace{-0.4cm}\sum_{\substack{1 \leq k \leq d \\ 1 \leq j \leq 2N_0+2}}h^{1+\frac{j}{2}}\overline{A_{l,k}^{(j)}}\phi_k\Big)_{L^2}=
\mathcal{O}(h^{N_0+\frac{5}{2}}\|u\|_{L^2}).$$
We get from (\ref{sl18}), (\ref{sl21}) and (\ref{sl75}) that the term factoring $h$ in
$$\sum_{0 \leq j,k \leq 2N_0+2}h^{1+\frac{j+k}{2}}\overline{a}_k^w(x,D_x)\psi_{j,l}^-+\sum_{\substack{1 \leq k \leq d \\ 1 \leq j \leq 2N_0+2}}h^{1+\frac{j}{2}}\overline{A_{l,k}^{(j)}}\phi_k,$$
is zero. It will therefore suffice to choose the Schwartz functions $\psi_{j,l}^-$, for $1 \leq l \leq d$, $1 \leq j \leq 2N_0+2$, such that 
$$\sum_{\substack{0 \leq k_1,k_2 \leq 2N_0+2 \\ k_1+k_2=j}}\overline{a}_{k_1}^w(x,D_x)\psi_{k_2,l}^-+\sum_{1 \leq k \leq d}\overline{A_{l,k}^{(j)}}\phi_k=0,$$
that is
\begin{equation}\label{eq7}
Q^*\psi_{j,l}^-=-\sum_{1 \leq k \leq d}\overline{A_{l,k}^{(j)}}\phi_k-\sum_{\substack{0 \leq k_1,k_2 \leq 2N_0+2 \\ k_1+k_2=j, \ k_1 \geq 1}}\overline{a}_{k_1}^w(x,D_x)\psi_{k_2,l}^-,
\end{equation}
for all $1 \leq l \leq d$, $1 \leq j \leq 2N_0+2$.
Assuming that the Schwartz functions $\psi_{k,l}^-$ have already been determined for all $0 \leq k \leq j-1$, $1 \leq l \leq d$, 
by using (\ref{sl80}) and (\ref{sl23bis}) we define the functions $(\psi_{j,l}^-)_{1 \leq l \leq d}$ as
\begin{align}\notag
\psi_{j,l}^-= & \ -(Q^*|_{V_2^{\perp}})^{-1}\pi_1\Big(\sum_{\substack{0 \leq k_1,k_2 \leq 2N_0+2 \\ k_1+k_2=j, \ k_1 \geq 1}}\overline{a}_{k_1}^w(x,D_x)\psi_{k_2,l}^-\Big)\\
= & \ -S^*\Big(\sum_{\substack{0 \leq k_1,k_2 \leq 2N_0+2 \\ k_1+k_2=j, \ k_1 \geq 1}}\overline{a}_{k_1}^w(x,D_x)\psi_{k_2,l}^-\Big). \label{sl82}
\end{align}
The next lemma establishes the identity 
$$\overline{A_{l,k}^{(j)}}=-\sum_{\substack{0 \leq k_1,k_2 \leq 2N_0+2 \\ k_1+k_2=j, \ k_1 \geq 1}}\big(\overline{a}_{k_1}^w(x,D_x)\psi_{k_2,l}^-,\phi_k\big)_{L^2},$$
which yields that equations (\ref{eq7}) are satisfied and therefore that equation (\ref{eq2}) holds.

\bigskip

\begin{lemma}\label{lem1.0}
The functions $\phi_{k,l}^+$, $\psi_{k,l}^-$ constructed above satisfy the identities
$$\sum_{\substack{0 \leq k_1,k_2 \leq 2N_0+2 \\ k_1+k_2=j, \ k_1 \geq 1}} \big(a_{k_1}^w(x,D_x)\phi_{k_2,l}^+,\psi_k\big)_{L^2}=\sum_{\substack{0 \leq k_1,k_2 \leq 2N_0+2 \\ k_1+k_2=j, \ k_1 \geq 1}}\big(\phi_l,\overline{a}_{k_1}^w(x,D_x)\psi_{k_2,k}^-\big)_{L^2},$$
for every $1 \leq j \leq 2N_0+2$, $1 \leq k,l \leq d$. Furthermore, the entries of the matrices $A_j=(A_{k,l}^{(j)})_{1 \leq k,l \leq d}$ are given by
$$A_{k,l}^{(j)}=\sum_{i=1}^j(-1)^{i}\sum_{\substack{1 \leq k_p \leq 2N_0+2 \\ k_1+...+k_i=j}} (a_{k_1}^wS a_{k_2}^wS\ ...\ a_{k_{i-1}}^wS a_{k_i}^w\phi_{l},\psi_k)_{L^2},$$
for all $1 \leq j \leq 2N_0+2$, $1 \leq k,l \leq d$.
\end{lemma}

\bigskip

\noindent
\begin{proof}
For $1 \leq k,l \leq d$, we have from (\ref{sl75}) and (\ref{eq6.5}) that
\begin{multline*}
\sum_{\substack{0 \leq k_1,k_2 \leq 2N_0+2 \\ k_1+k_2=j, \ k_1 \geq 1}} (a_{k_1}^w\phi_{k_2,l}^+,\psi_k)_{L^2}=\sum_{\substack{0 \leq k_1,k_2 \leq 2N_0+2 \\ k_1+k_2=j, \ k_1 \geq 1}} (\phi_{k_2,l}^+,\overline{a}_{k_1}^w\psi_k)_{L^2}\\ 
=(\phi_{l},\overline{a}_{j}^w\psi_k)_{L^2} -\sum_{\substack{1 \leq k_1,k_2 \leq 2N_0+2 \\ k_1+k_2=j}} \sum_{\substack{0 \leq k_3,k_4 \leq 2N_0+2 \\ k_3+k_4=k_2, \ k_3 \geq 1}}(S a_{k_3}^w\phi_{k_4,l}^+,\overline{a}_{k_1}^w\psi_k)_{L^2}.
\end{multline*}
We may write 
\begin{multline*}
\sum_{\substack{1 \leq k_1,k_2 \leq 2N_0+2 \\ k_1+k_2=j}} \sum_{\substack{0 \leq k_3,k_4 \leq 2N_0+2 \\ k_3+k_4=k_2, \ k_3 \geq 1}}(S a_{k_3}^w\phi_{k_4,l}^+,\overline{a}_{k_1}^w\psi_k)_{L^2}\\
=\sum_{\substack{1 \leq k_1,k_2 \leq 2N_0+2 \\ k_1+k_2=j}} (\phi_{l},\overline{a}_{k_2}^wS^*\overline{a}_{k_1}^w\psi_k)_{L^2}
+\sum_{\substack{1 \leq k_1,k_3,k_4 \leq 2N_0+2 \\ k_1+k_3+k_4=j}} (\phi_{k_4,l}^+,\overline{a}_{k_3}^wS^*\overline{a}_{k_1}^w\psi_k)_{L^2},
\end{multline*}
for all $1 \leq j \leq 2N_0+2$, which yields
\begin{multline*}
\sum_{\substack{0 \leq k_1,k_2 \leq 2N_0+2 \\ k_1+k_2=j, \ k_1 \geq 1}} (a_{k_1}^w\phi_{k_2,l}^+,\psi_k)_{L^2} =(\phi_{l},\overline{a}_{j}^w\psi_k)_{L^2} -
\sum_{\substack{1 \leq k_1,k_2 \leq 2N_0+2 \\ k_1+k_2=j}} (\phi_{l},\overline{a}_{k_1}^wS^*\overline{a}_{k_2}^w\psi_k)_{L^2}\\ -\sum_{\substack{1 \leq k_1,k_2,k_3 \leq 2N_0+2 \\ k_1+k_2+k_3=j}} (\phi_{k_1,l}^+,\overline{a}_{k_2}^wS^*\overline{a}_{k_3}^w\psi_k)_{L^2}.
\end{multline*}
By using the definition (\ref{eq6.5}) of the functions $\phi_{k,l}^+$ and iterating this process, we obtain 
\begin{multline*}
\sum_{\substack{0 \leq k_1,k_2 \leq 2N_0+2 \\ k_1+k_2=j, \ k_1 \geq 1}} (a_{k_1}^w\phi_{k_2,l}^+,\psi_k)_{L^2} \\ = 
\sum_{i=1}^j(-1)^{i+1}\sum_{\substack{1 \leq k_p \leq 2N_0+2 \\ k_1+...+k_i=j}} (\phi_{l},\overline{a}_{k_1}^wS^*\overline{a}_{k_2}^wS^*...\ \overline{a}_{k_{i-1}}^wS^*\overline{a}_{k_i}^w\psi_k)_{L^2}.
\end{multline*}
On the other hand, from (\ref{sl75}) and (\ref{sl82}) it follows that
\begin{multline*}
\sum_{\substack{0 \leq k_1,k_2 \leq 2N_0+2 \\ k_1+k_2=j, \ k_1 \geq 1}}(\phi_l,\overline{a}_{k_1}^w\psi_{k_2,k}^-)_{L^2}=\sum_{\substack{0 \leq k_1,k_2 \leq 2N_0+2 \\ k_1+k_2=j, \ k_1 \geq 1}}(a_{k_1}^w\phi_l,\psi_{k_2,k}^-)_{L^2}\\
=(a_{j}^w\phi_l,\psi_{k})_{L^2}-\sum_{\substack{1 \leq k_1,k_2 \leq 2N_0+2 \\ k_1+k_2=j}}\sum_{\substack{0 \leq k_3,k_4 \leq 2N_0+2 \\ k_3+k_4=k_2, \ k_3 \geq 1}}(a_{k_1}^w\phi_l,S^*\overline{a}_{k_3}^w\psi_{k_4,k}^-)_{L^2}.
\end{multline*}
Since we may write 
\begin{multline*}
\sum_{\substack{1 \leq k_1,k_2 \leq 2N_0+2 \\ k_1+k_2=j}}\sum_{\substack{0 \leq k_3,k_4 \leq 2N_0+2 \\ k_3+k_4=k_2, \ k_3 \geq 1}}(a_{k_1}^w\phi_l,S^*\overline{a}_{k_3}^w\psi_{k_4,k}^-)_{L^2}\\ 
=\sum_{\substack{1 \leq k_1,k_2 \leq 2N_0+2 \\ k_1+k_2=j}}(a_{k_2}^wS a_{k_1}^w\phi_l,\psi_{k})_{L^2}+\sum_{\substack{1 \leq k_1,k_3,k_4 \leq 2N_0+2 \\ k_1+k_3+k_4=j}}(a_{k_3}^wS a_{k_1}^w\phi_l,\psi_{k_4,k}^-)_{L^2},
\end{multline*}
we get
\begin{multline*}
\sum_{\substack{0 \leq k_1,k_2 \leq 2N_0+2 \\ k_1+k_2=j, \ k_1 \geq 1}}(\phi_l,\overline{a}_{k_1}^w\psi_{k_2,k}^-)_{L^2}=(a_{j}^w\phi_l,\psi_{k})_{L^2} -\sum_{\substack{1 \leq k_1,k_2 \leq 2N_0+2 \\ 
k_1+k_2=j}}(a_{k_1}^wS a_{k_2}^w\phi_l,\psi_{k})_{L^2}\\ -\sum_{\substack{1 \leq k_1,k_2,k_3 \leq 2N_0+2 \\ k_1+k_2+k_3=j}}(a_{k_1}^wS a_{k_2}^w\phi_l,\psi_{k_3,k}^-)_{L^2}.
\end{multline*}
By using the definition (\ref{sl82}) of the functions $\psi_{k,l}^-$ and iterating this process, we obtain that 
\begin{multline*}
\sum_{\substack{0 \leq k_1,k_2 \leq 2N_0+2 \\ k_1+k_2=j, \ k_1 \geq 1}}(\phi_l,\overline{a}_{k_1}^w\psi_{k_2,k}^-)_{L^2}\\
=\sum_{i=1}^j(-1)^{i+1}\sum_{\substack{1 \leq k_p \leq 2N_0+2 \\ k_1+...+k_i=j}} (a_{k_1}^wS a_{k_2}^wS\ ...\ a_{k_{i-1}}^wS a_{k_i}^w\phi_{l},\psi_k)_{L^2}.
\end{multline*}
As
$$(\phi_{l},\overline{a}_{k_1}^wS^*\overline{a}_{k_2}^wS^*...\ \overline{a}_{k_{i-1}}^wS^*\overline{a}_{k_i}^w\psi_k)_{L^2}=
(a_{k_i}^wS a_{k_{i-1}}^wS\ ...\ a_{k_{2}}^wS a_{k_1}^w\phi_{l},\psi_k)_{L^2},
$$
we hence conclude from (\ref{eq6}) that 
$$A_{k,l}^{(j)}=-\sum_{\substack{0 \leq k_1,k_2 \leq 2N_0+2 \\ k_1+k_2=j, \ k_1 \geq 1}} (a_{k_1}^w\phi_{k_2,l}^+,\psi_k)_{L^2}=-
\sum_{\substack{0 \leq k_1,k_2 \leq 2N_0+2 \\ k_1+k_2=j, \ k_1 \geq 1}}(\phi_l,\overline{a}_{k_1}^w\psi_{k_2,k}^-)_{L^2},$$
for all $1 \leq j \leq 2N_0+2$.
\end{proof}

\noindent
Writing $u=u_1+u_2 \in L^2(\rr^n)$, with $(u_1,u_2) \in V_1 \times V_1^{\perp}$, we finally obtain 
from (\ref{sl17}), (\ref{sl24}) and (\ref{sl75}) that equation (\ref{eq1}) readily holds:
\begin{multline*}
S Qu+\sum_{k=1}^d(u,\phi_k)_{L^2}\left(\sum_{j=0}^{2N_0+2}\phi_{j,k}^+h^{\frac{j}{2}}\right)=u_2+\sum_{k=1}^d(u,\phi_k)_{L^2}\phi_{k}+\mathcal{O}(h^{1/2}\|u\|_{L^2})\\
=u_2+u_1+\mathcal{O}(h^{1/2}\|u\|_{L^2})=u+\mathcal{O}(h^{1/2}\|u\|_{L^2}).
\end{multline*}
We shall now use the Grushin-reduction (\ref{eq4}), (\ref{eq3}), (\ref{eq2}), (\ref{eq1}) in order to prove Theorem~\ref{th1}.

Let $\Omega$ be a compact subset of $K^{2N_0+2}$.
We first assume that there exist $c_0>0$, $0<h_0\leq 1$ such that for all $u \in L^2(\rr^n)$, $0 <h \leq h_0$, and $(z_1,...,z_{2N_0+2}) \in \Omega$,
\begin{equation}\label{eq8}
 \|Pu-hz(h)u\|_{L^2} \geq c_0 h^{\frac{N_0}{2}+1}\|u\|_{L^2}.
\end{equation}
From (\ref{eq10.1}) and (\ref{eq8}) we get that for any given $u_- \in \cc^d$,
\begin{multline}\label{sl100}
 c_0 h^{\frac{N_0}{2}+1}\|E_+u_-\|_{L^2} \leq \|(P-hz(h))E_+u_-\|_{L^2} \\ 
\leq  \Big\|\sum_{k=0}^{2N_0+2}h^{1+\frac{k}{2}} a_k^w(x,D_x)E_+u_-+R_{-}E_{\pm}u_-\Big\|_{L^2} +\|R_{-}E_{\pm}u_-\|_{L^2}\\  +   
h^{\frac{N_0}{2}+\frac{3}{2}}\Big\|\sum_{|\alpha|=N_0+3}\textrm{Op}^w(X^{\alpha}R_{\alpha}(X;h))E_+u_-\Big\|_{L^2}+\mathcal{O}(h^{[\frac{N_0}{2}]+2})\|E_+u_-\|_{L^2}. 
\end{multline}
By using a little symbolic calculus in the Weyl quantization, we readily obtain from (\ref{sl72}) and the exact formula (\cite{hormander}, Theorem~18.5.4), 
\begin{equation}\label{kk1}
X^{\alpha} \#^w R_{\alpha}=\sum_{p=0}^{|\alpha|}\frac{1}{p!}\Big(\frac{1}{2i}\sigma(\partial_{X_1},\partial_{X_2})\Big)^pX_1^{\alpha}R_{\alpha}(X_2;h)\Big{|}_{X_1=X_2=X},
\end{equation}
that the operator $\textrm{Op}^w(X^{\alpha}R_{\alpha}(X;h))$ may be written as
\begin{equation}\label{sl90}
\textrm{Op}^w(X^{\alpha}R_{\alpha}(X;h))=\sum_{\beta \leq \alpha}\textrm{Op}^w(\tilde{R}_{\beta}(X;h))\textrm{Op}^w(X^{\beta}),
\end{equation}
for some symbols $\tilde{R}_{\beta}$ belonging to the class $S(1)$.
It follows from (\ref{eq10.2}) and (\ref{sl90}) that 
\begin{align}\label{sl101}
& \ \|\textrm{Op}^w(X^{\alpha}R_{\alpha}(X;h))E_+u_-\|_{L^2} \\ \notag & \qquad \leq \sum_{\substack{1 \leq k \leq d \\ 0 \leq j \leq 2N_0+2}}|u_-(k)|h^{\frac{j}{2}}\| \textrm{Op}^w(X^{\alpha}R_{\alpha}(X;h))\phi_{j,k}^+\|_{L^2}\\ \notag & \qquad
\lesssim |u_-| \sum_{\substack{1 \leq k \leq d \\ 0 \leq j \leq 2N_0+2\\ \beta \leq \alpha}}\|\textrm{Op}^w(X^{\beta})\phi_{j,k}^+\|_{L^2}=\mathcal{O}(1)|u_-|,
\end{align}
since the functions $\phi_{j,k}^+$ belong to $\mathscr{S}(\rr^n)$. 
As $\|R_-\|_{\mathscr{L}(\cc^d,L^2)}=\mathcal{O}(1)$, it follows from (\ref{eq3}), (\ref{sl100}) and (\ref{sl101}) that 
\begin{equation}\label{sl102}
h^{\frac{N_0}{2}+1}\|E_+u_-\|_{L^2} \lesssim  |E_{\pm}u_-|+ \mathcal{O}(h^{\frac{N_0}{2}+\frac{3}{2}})|u_-|.
\end{equation}   
Since $\|R_+\|_{\mathscr{L}(L^2,\cc^d)}=\mathcal{O}(1)$, we get from (\ref{eq4}) that 
\begin{equation}\label{sl103}
|u_-| \leq |R_+E_+u_-|+\mathcal{O}(h^{1/2})|u_-| \leq \|E_+u_-\|_{L^2}+\mathcal{O}(h^{1/2})|u_-|,
\end{equation}
whence from (\ref{sl102}) and (\ref{sl103}) we obtain 
$$h^{\frac{N_0}{2}+1}|u_-| \lesssim  |E_{\pm}u_-|+ \mathcal{O}(h^{\frac{N_0}{2}+\frac{3}{2}})|u_-|,$$
that is, there exist constants $c_0>0$, $0<h_0\leq 1$, such that
\begin{multline}\label{eq10.4}
\forall u_- \in \cc^d, \forall \ 0 <h \leq h_0, \forall (z_1,...,z_{2N_0+2}) \in \Omega,    |E_{\pm}u_-| \geq c_0 h^{\frac{N_0}{2}+1}|u_-|.
\end{multline}
This ends the proof of the first implication. 

We shall now prove the converse implication. We therefore assume that the estimate (\ref{eq10.4}) holds. 
It follows from (\ref{sl5}),  (\ref{sl9}), (\ref{sl11}), (\ref{sl13}) and (\ref{sl16}) that
\begin{align}
p(h^{1/2}X;h)-hz(h)= & \ p_0(h^{1/2}X)+hp_1(h^{1/2}X)-hz_0 \ \textrm{mod } S(h^{3/2}) \notag \\
=& \ h\big(q(X)+p_1(0)-z_0\big)+r_{0,h}(X)+r_{1,h}(X) \ \textrm{mod } S(h^{3/2}) \notag \\
=& \ ha_0(X)+r_{0,h}(X)+r_{1,h}(X) \ \textrm{mod } S(h^{3/2}),\label{eq10.5}
\end{align}
with 
\begin{equation}\label{k2.1}
r_{0,h}(X)=\sum_{|\alpha|=3}\frac{3}{\alpha!}X^{\alpha}h^{3/2}\int_0^1(1-t)^{2}p_0^{(\alpha)}(th^{1/2}X)dt,
\end{equation} 
\begin{equation}\label{k2.2}
r_{1,h}(X)=\sum_{|\alpha|=1}X^{\alpha}h^{3/2}\int_0^1p_1^{(\alpha)}(th^{1/2}X)dt.
\end{equation}
Let $\chi_0 \in C_0^{\infty}(\rr^{2n})$ be a cutoff function satisfying $0\leq\chi_0\leq 1$ and 
\begin{equation}\label{kk4}
\textrm{supp } \chi_0 \subset \{X \in \rr^{2n} ; |X| \leq 2\}, \quad \chi_0=1 \textrm{ on } \{X \in \rr^{2n} ; |X| \leq 1\},
\end{equation}
and let $A \gg 1$ be a large positive constant to be chosen later on.
Setting
\begin{equation}\label{k1.1}
M_0=\chi_0^w(Ah^{1/2}x,Ah^{1/2}D_x),
\end{equation}
it follows from (\ref{sl17}) and (\ref{eq1}) that for all $u \in \mathscr{S}(\rr^n)$,
\begin{multline*}
h\|u\|_{L^2} \leq h\|E_+R_+u\|_{L^2}+h\|SQu\|_{L^2} +\mathcal{O}(h^{3/2})\|u\|_{L^2}
\lesssim  h|R_+u|\\ + h\|SQM_0u\|_{L^2} + h\|SQ(1-M_0)u\|_{L^2} 
+\mathcal{O}(h^{3/2})\|u\|_{L^2},
\end{multline*}
since $\|E_+\|_{\mathscr{L}(\cc^d,L^2)}=\mathcal{O}(1)$. From (\ref{sl24}) one has 
$$\|SQ(1-M_0)u\|_{L^2}=\|\pi_1(1-M_0)u\|_{L^2} \leq \|(1-M_0)u\|_{L^2},$$
whence 
\begin{equation}\label{eq12.2}
h\|u\|_{L^2} \lesssim  h|R_+u|+h\|SQM_0u\|_{L^2} + h\|(1-M_0)u\|_{L^2}.
\end{equation}
Observing that $\|S\|_{\mathscr{L}(L^2)}=\mathcal{O}(1)$ and $\|M_0\|_{\mathscr{L}(L^2)}=\mathcal{O}(1)$, when $0<h \leq A^{-2} \leq 1$,
we then deduce from (\ref{sl17}) and (\ref{eq10.5}) that 
\begin{multline}\label{k1.2}
h\|SQM_0u\|_{L^2} \leq  \|S(P-hz(h))M_0u\|_{L^2}+\|S r_{0,h}^wM_0u\|_{L^2}\\ +\|S r_{1,h}^wM_0u\|_{L^2}+\mathcal{O}(h^{3/2})\|u\|_{L^2},
\end{multline}
that in turn yields 
\begin{multline}\label{k1.3}
\|S(P-hz(h))M_0u\|_{L^2} \lesssim \|(P-hz(h))M_0u\big\|_{L^2} \\ \lesssim  \|M_0(P-hz(h))u\|_{L^2} +\|[P,M_0]u\|_{L^2} \lesssim  \|Pu-hz(h)u\|_{L^2}+\|[P,M_0]u\|_{L^2},
\end{multline}
which, along with (\ref{k1.2}), gives
\begin{multline}\label{eq12.1}
h\|SQM_0u\|_{L^2} \lesssim  \|Pu-hz(h)u\|_{L^2}+\|[P,M_0]u\|_{L^2}\\  +\|S r_{0,h}^wM_0u\|_{L^2}+\|S r_{1,h}^wM_0u\|_{L^2}+\mathcal{O}(h^{3/2})\|u\|_{L^2}.
\end{multline}
We shall need the following technical lemma.

\bigskip

\begin{lemma}\label{lem1.1}
We have 
$$\|S r_{1,h}^wM_0u\|_{L^2}=\|S r_{1,h}^w\chi_0^w(Ah^{1/2}X)u\|_{L^2}=\mathcal{O}\Big(\frac{h}{A}\Big)\|u\|_{L^2}+\mathcal{O}_A(h^2)\|u\|_{L^2},$$
when $0<h \leq A^{-2} \leq 1$.
\end{lemma}

\bigskip

\begin{proof} 
Since $\|S\|_{\mathscr{L}(L^2)}=\mathcal{O}(1)$, we notice in the first place that 
$$\|S r_{1,h}^wM_0u\|_{L^2} \leq \|r_{1,h}^wM_0u\|_{L^2}.$$
By referring to Section~\ref{setting} for the definitions of the symbol classes, we then show that if $R_1(\cdot;h) \in S(1)$, $R_2(\cdot;h) \in S(\langle X \rangle^{-1})$, then there exists $R_3(\cdot;h) \in S(h)$ such that 
\begin{multline}\label{kk2}
\Bigl(h^{1/2}X_jR_1(h^{1/2}X;h)\Bigr) \#^w R_2(h^{1/2}X;h)\\ =h^{1/2}X_jR_1(h^{1/2}X;h)R_2(h^{1/2}X;h)+R_3(h^{1/2}X;h).
\end{multline}
Indeed, we deduce from (\ref{kk1}) that 
\begin{align*}
& \ \Bigl(h^{1/2}X_jR_1(h^{1/2}X;h)\Bigr) \#^w R_2(h^{1/2}X;h)\\ 
=& \ R_1(h^{1/2}X;h)\#^w (h^{1/2}X_j) \#^w R_2(h^{1/2}X;h)+\frac{h}{2i}R_4(h^{1/2}X;h)\\ 
= & \ R_1(h^{1/2}X;h)\#^w \Bigl(h^{1/2}X_jR_2(h^{1/2}X;h)\Bigr)+\frac{h}{2i}R_4(h^{1/2}X;h)+\frac{h}{2i}R_5(h^{1/2}X;h),
\end{align*}
with
$$R_4(h^{1/2}X;h)=\Bigl(\{X_j,R_1(\cdot;h)\}(h^{1/2}X)\Bigr)\#^w R_2(h^{1/2}X;h),$$
$$R_5(h^{1/2}X;h)=R_1(h^{1/2}X;h) \#^w \Bigl(\{X_j,R_2(\cdot;h)\}(h^{1/2}X)\Bigr).$$
By using the symbolic calculus (\cite{DS}, Chapter~7), we notice that $R_4(\cdot;h) \in S(\langle X\rangle^{-1})$, because $\{X_j,R_1(\cdot;h)\} \in S(1)$ and $R_2(\cdot;h) \in S(\langle X \rangle^{-1})$. We also notice that 
$R_5(\cdot;h) \in S(\langle X\rangle^{-1})$, because $R_1(\cdot;h) \in S(1)$ and $\{X_j,R_2(\cdot;h)\} \in S(\langle X\rangle^{-1})$. It follows that 
\begin{multline*}
\Bigl(h^{1/2}X_jR_1(h^{1/2}X;h)\Bigr) \#^w R_2(h^{1/2}X;h) \\
=R_1(h^{1/2}X;h)\#^w \Bigl(h^{1/2}X_jR_2(h^{1/2}X;h)\Bigr)+R_6(h^{1/2}X;h),
\end{multline*}
with $R_6(\cdot;h) \in S(h)$. Next, since $X_jR_2(X;h) \in S(1)$, another use of the symbolic calculus in the class $S(1)$ gives that 
\begin{multline*}
R_1(h^{1/2}X;h)\#^w \Bigl(h^{1/2}X_jR_2(h^{1/2}X;h)\Bigr)\\
=h^{1/2}X_jR_1(h^{1/2}X;h)R_2(h^{1/2}X;h)+R_7(h^{1/2}X;h),
\end{multline*}
with $R_7(\cdot;h) \in S(h)$, and this proves (\ref{kk2}). 

We next notice from (\ref{sl5}) and (\ref{kk4}) that
\begin{equation}\label{bb1}
\int_0^1p_1^{(\alpha)}(tX)dt \in S(1), \quad \chi_0(AX) \in S(\mathcal{O}_A(\langle X \rangle^{-1})).
\end{equation}
It therefore follows from (\ref{k2.2}), (\ref{kk2}) and (\ref{bb1}) that 
\begin{align*}
 \frac{A}{h} r_{1,h} &\#^w \chi_0(Ah^{1/2}X)= \Big(\sum_{|\alpha|=1}(Ah^{1/2}X)^{\alpha}\int_0^1p_1^{(\alpha)}(th^{1/2}X)dt\Big) \#^w \chi_0(Ah^{1/2}X)\\
=& \ \sum_{|\alpha|=1}(Ah^{1/2}X)^{\alpha}\chi_0(Ah^{1/2}X)\int_0^1p_1^{(\alpha)}(th^{1/2}X)dt+R_8(h^{1/2}X;h,A),
\end{align*}
with $R_8(\cdot;h,A) \in S(\mathcal{O}_A(h))$.
Since the symbol 
$$\sum_{|\alpha|=1}(Ah^{1/2}X)^{\alpha}\chi_0(Ah^{1/2}X)\int_0^1p_1^{(\alpha)}(th^{1/2}X)dt,$$
belongs to the class $S(1)$ uniformly with respect to the parameters when $0<h \leq A^{-2} \leq 1$, we thus get
$$\|r_{1,h}^wM_0u\|_{L^2}=\|r_{1,h}^w\chi_0^w(Ah^{1/2}X)u\|_{L^2}=\mathcal{O}\Big(\frac{h}{A}\Big)\|u\|_{L^2}+\mathcal{O}_A(h^2)\|u\|_{L^2},$$
 when $0<h \leq A^{-2} \leq 1$, and this concludes the proof of the lemma.
\end{proof}

\noindent
We shall also need the following technical result.

\bigskip

\begin{lemma}\label{lem1.2}
We have 
$$\|S r_{0,h}^wM_0u\|_{L^2}=\mathcal{O}\Big(\frac{h}{A}\Big)\|u\|_{L^2}+\mathcal{O}_A(h^2)\|u\|_{L^2},$$
when $0<h \leq A^{-2} \leq 1$.
\end{lemma}

\medskip

\begin{proof}
From (\ref{k2.1}) we have 
\begin{multline}\label{dd1}
\frac{A}{h} r_{0,h}(X) \#^w \chi_0(Ah^{1/2}X) \\ =\Big(\sum_{|\alpha|=3}\frac{3}{\alpha!}X^{\alpha}Ah^{1/2}\int_0^1(1-t)^{2}p_0^{(\alpha)}(th^{1/2}X)dt\Big) \#^w \chi_0(Ah^{1/2}X).
\end{multline}
Observing from (\ref{kk1}) that if $R(\cdot;h) \in S(1)$ then for each $\alpha\in\nn^n$ 
there exist symbols $R_{\beta}(\cdot;h) \in S(1)$, with $\beta \leq \alpha$, $|\beta|<|\alpha|$, $\beta \in \nn^n$, such that 
\begin{equation}\label{bb3}
X^{\alpha}R(h^{1/2}X;h)= X^{\alpha}\#^w R(h^{1/2}X;h)+\sum_{\substack{\beta \leq \alpha \\ |\beta|<|\alpha|}}h^{\frac{|\alpha|-|\beta|}{2}}X^{\beta}R_{\beta}(h^{1/2}X;h),
\end{equation}
by induction we readily have that  if $R(\cdot;h) \in S(1)$, the for each $\alpha\in\nn^n$ there exist symbols $R_{\beta}(\cdot;h) \in S(1)$, with 
$\beta \leq \alpha$, $|\beta|<|\alpha|$, $\beta \in \nn^n$, such that 
\begin{multline}\label{bb4}
X^{\alpha}R(h^{1/2}X;h)= X^{\alpha}\#^w R(h^{1/2}X;h)\\ +\sum_{\substack{\beta \leq \alpha \\ |\beta|<|\alpha|}}h^{\frac{|\alpha|-|\beta|}{2}}X^{\beta}\#^w R_{\beta}(h^{1/2}X;h).
\end{multline}
We deduce from (\ref{sl5}), (\ref{dd1}) and (\ref{bb4}) that there exist symbols $R_{\beta}(\cdot;h,A)\in S(\mathcal{O}_A(1))$, for $|\beta| \leq 2$, such that
\begin{align*}
& \ \frac{A}{h} r_{0,h}(X) \#^w \chi_0(Ah^{1/2}X)\\
=& \ \sum_{|\alpha|=3}\frac{3}{\alpha!}X^{\alpha}Ah^{1/2}\#^w \Big(\int_0^1(1-t)^{2}p_0^{(\alpha)}(th^{1/2}X)dt\Big) \#^w \chi_0(Ah^{1/2}X)\\
& \  + h\sum_{|\beta| \leq 2}h^{\frac{2-|\beta|}{2}}X^{\beta}\#^w R_{\beta}(h^{1/2}X;h,A)\#^w \chi_0(Ah^{1/2}X).
\end{align*}
The symbolic calculus shows that there exists $r_1(\cdot;h,A) \in S(\mathcal{O}_A(\langle X \rangle^2))$ such that
\begin{equation}\label{bb5} 
r_1(X;h,A)=\sum_{|\beta| \leq 2}h^{\frac{2-|\beta|}{2}}X^{\beta}\#^w R_{\beta}(h^{1/2}X;h,A)\#^w \chi_0(Ah^{1/2}X),
\end{equation}
since we have $X^{\beta} \in S(\langle X \rangle^2)$ when $|\beta| \leq 2$, $R_{\beta}(h^{1/2}X;h,A) \in S(\mathcal{O}_A(1))$, and 
$\chi_0(Ah^{1/2}X) \in S(1)$ uniformly with respect to the parameters $h, A$, when $0<h \leq A^{-2} \leq 1$.
On the other hand, we get from (\ref{sl5}), (\ref{kk4}) and another use of the symbolic calculus 
that there exists a symbol $r_2(\cdot;h,A) \in S(\mathcal{O}_A(\langle X \rangle^{-\infty}))$ such that
\begin{multline}\label{bb6}
\Big(\int_0^1(1-t)^{2}p_0^{(\alpha)}(th^{1/2}X)dt\Big) \#^w \chi_0(Ah^{1/2}X)\\ =\chi_0(Ah^{1/2}X)\int_0^1(1-t)^{2}p_0^{(\alpha)}(th^{1/2}X)dt+hr_2(h^{1/2}X;h,A).
\end{multline}
It follows that 
\begin{multline}\label{bb7}
\frac{A}{h} r_{0,h}(X) \#^w \chi_0(Ah^{1/2}X)=h\sum_{|\alpha|=3}\frac{3}{\alpha!}X^{\alpha}Ah^{1/2}\#^w r_2(h^{1/2}X;h,A)\\
+\sum_{|\alpha|=3}\frac{3}{\alpha!}X^{\alpha}Ah^{1/2}\#^w \Big(\chi_0(Ah^{1/2}X)\int_0^1(1-t)^{2}p_0^{(\alpha)}(th^{1/2}X)dt\Big)+hr_1(X;h,A).
\end{multline}
Using (\ref{kk1}) shows that there exist symbols $r_3(\cdot;h,A)$, $r_4(\cdot;h,A) \in S(\mathcal{O}_A(\langle X \rangle^{-\infty}))$, 
and $r_5(\cdot;h,A) \in S(\mathcal{O}_A(\langle X \rangle^2))$ such that
\begin{align*}
& \ \Bigl(h^{1/2}X_{j_1}X_{j_2}X_{j_3}\Bigr)\#^w r_2(h^{1/2}X;h,A)= \Bigl(X_{j_1}X_{j_2}\Bigr) \#^w \Bigl(h^{1/2}X_{j_3}\Bigr)\#^w r_2(h^{1/2}X;h,A)\\
& \ +\frac{i}{2}\Bigl(h^{1/2}\{X_{j_1}X_{j_2},X_{j_3}\}\Bigr)\#^w r_2(h^{1/2}X;h,A) =\Bigl((X_{j_1}X_{j_2})\#^w r_3(h^{1/2}X;h,A)\Bigr)\\
& \ +r_4(h^{1/2}X;h,A)=r_5(X;h,A),
\end{align*}
when $0<h \leq A^{-2} \leq 1$.
It follows from (\ref{bb7}) and the previous identity that there exists a symbol $r_6(\cdot;h,A) \in S(\mathcal{O}_A(\langle X \rangle^2))$ such that
\begin{multline}\label{bb8}
\frac{A}{h} r_{0,h}(X) \#^w \chi_0(Ah^{1/2}X)\\
=\sum_{|\alpha|=3}\frac{3}{\alpha!}X^{\alpha}Ah^{1/2}\#^w \Big(\chi_0(Ah^{1/2}X)\int_0^1(1-t)^{2}p_0^{(\alpha)}(th^{1/2}X)dt\Big)+hr_6(X;h,A),
\end{multline}
when $0<h \leq A^{-2} \leq 1$.
When $|\alpha|=3$, formula (\ref{kk1}) once more gives that there exist symbols $R_{\beta}(\cdot;h,A) \in S(\mathcal{O}_A(1))$, with $|\beta| \leq 2$, such that
\begin{align*}
& \ X^{\alpha}\#^w \Big(\chi_0(Ah^{1/2}X)\int_0^1(1-t)^{2}p_0^{(\alpha)}(th^{1/2}X)dt\Big)\\
= & \  X^{\alpha}\chi_0(Ah^{1/2}X)\int_0^1(1-t)^{2}p_0^{(\alpha)}(th^{1/2}X)dt\\
& \quad +h^{1/2}\sum_{|\beta| \leq 2}h^{\frac{2-|\beta|}{2}}X^{\beta}R_{\beta}(h^{1/2}X;h,A),
\end{align*}
because 
$$\chi_0(Ah^{1/2}X)\int_0^1(1-t)^{2}p_0^{(\alpha)}(th^{1/2}X)dt \in S(\mathcal{O}_A(1)),$$
when $0<h \leq A^{-2} \leq 1$.
It follows from (\ref{bb8}) and the previous identity that there exists a symbol $r_7(\cdot;h,A) \in S(\mathcal{O}_A(\langle X \rangle^2))$ such that 
\begin{multline}\label{bb9}
\frac{A}{h} r_{0,h}(X) \#^w \chi_0(Ah^{1/2}X)\\
=\sum_{|\alpha|=3}\frac{3}{\alpha!}X^{\alpha}Ah^{1/2}\chi_0(Ah^{1/2}X)\int_0^1(1-t)^{2}p_0^{(\alpha)}(th^{1/2}X)dt+hr_7(X;h,A),
\end{multline}
because
$$h^{\frac{2-|\beta|}{2}}X^{\beta}R_{\beta}(h^{1/2}X;h,A) \in S(\mathcal{O}_A(\langle X \rangle^2)),$$
when $|\beta| \leq 2$ and $0<h \leq A^{-2} \leq 1$.
Consider then the symbol
\begin{equation}\label{bb10}
r_{2,h}(X)=\sum_{|\alpha|=3}\frac{3}{\alpha!}X^{\alpha}Ah^{1/2}\chi_0(Ah^{1/2}X)\int_0^1(1-t)^{2}p_0^{(\alpha)}(th^{1/2}X)dt,
\end{equation}
which may be written as
$$r_{2,h}(X)=\sum_{|\alpha_1|=2, |\alpha_2|=1} X^{\alpha_1}(Ah^{1/2}X)^{\alpha_2}\chi_0(Ah^{1/2}X)p_{\alpha_1,\alpha_2}(h^{1/2}X),$$
for some symbols $p_{\alpha_1,\alpha_2}$ belonging to the class $S(1)$, since $p_0 \in S(1)$.
We therefore deduce from (\ref{kk4}) that the symbol $r_{2,h}$ belongs to the class $S(\langle X \rangle^2)$ uniformly 
with respect to the parameters when $0<h \leq A^{-2} \leq 1$. 
By using the fact that the symbol of the operator $S$ belongs to the class $\mathbf{S}^{-2}$, we obtain from (\ref{bb9}) and (\ref{bb10}) that
$$S r_{0,h}^wM_0=S r_{0,h}^w\chi_0^w(Ah^{1/2}X)=\frac{h}{A}r_{3,h}^w+h^2r_{4,h}^w,$$ 
for some symbols $r_{3,h} \in S(1)$, $r_{4,h} \in S(\mathcal{O}_A(1))$ uniformly with respect to the parameters when $0<h \leq A^{-2} \leq 1$.
It follows that 
$$\|S r_{0,h}^wM_0u\|_{L^2}  \lesssim \frac{h}{A}\|u\|_{L^2}
+\mathcal{O}_A(h^2)\|u\|_{L^2},$$
when $0<h \leq A^{-2} \leq 1$. This ends the proof of Lemma~\ref{lem1.2}.
\end{proof}

\medskip

\noindent
We now resume the proof of Theorem~\ref{th1} and deduce from (\ref{eq12.1}), Lemmas~\ref{lem1.1} and~\ref{lem1.2} that 
\begin{multline}\label{bb12}
h\|SQM_0u\|_{L^2} \lesssim  \|Pu-hz(h)u\|_{L^2}+ \|[P,M_0]u\|_{L^2}\\  +\mathcal{O}\Big(\frac{h}{A}\Big)\|u\|_{L^2}+\mathcal{O}_A(h^2)\|u\|_{L^2}+\mathcal{O}(h^{3/2})\|u\|_{L^2},
\end{multline}
when $0<h \leq A^{-2} \leq 1$. Then, from (\ref{eq12.2}) and (\ref{bb12}), we get 
\begin{multline*}
 h\|u\|_{L^2} \lesssim   \|Pu-hz(h)u\|_{L^2}+ \|[P,M_0]u\|_{L^2}+  h|R_+u| \\ +h\|(1-M_0)u\|_{L^2}  + \mathcal{O}\Big(\frac{h}{A}\Big)\|u\|_{L^2} +
\mathcal{O}_A(h^2)\|u\|_{L^2}+\mathcal{O}(h^{3/2})\|u\|_{L^2},
\end{multline*}
when $0<h \leq A^{-2} \leq 1$. 
We next choose the large parameter $A \gg 1$ to control the term $\mathcal{O}(\frac{h}{A})\|u\|_{L^2}$ by the left-hand side of the previous estimate. 
With this definitive choice fixing the parameters $A_0 \geq 1$, $0<h_0\ll 1$, we hence obtain
\begin{multline}\label{eq12.4}
 h\|u\|_{L^2} \lesssim   \|Pu-hz(h)u\|_{L^2}+\|[P,M_0]u\|_{L^2}+  h|R_+u| \\ +h\|(1-M_0)u\|_{L^2}  +\mathcal{O}(h^{3/2})\|u\|_{L^2},
\end{multline}
when $0<h \leq h_0$.
Notice from (\ref{sl5}) and (\ref{kk4}) that the Weyl symbol of the operator 
$$[P,M_0]=\bigl[P,\chi_0^w(A_0h^{1/2}X)\bigr],$$ 
is given by
\begin{multline*}
\frac{1}{i}\big\{p(h^{1/2}X;h),\chi_0(A_0h^{1/2}X)\big\} \textrm{ mod } S(h^2)\\
=\frac{h}{i}\big\{p_0,\chi_0(A_0\cdot)\big\}(h^{1/2}X)  \textrm{ mod } S(h^2),
\end{multline*}
whence it follows that 
\begin{equation}\label{bb13}
\|[P,M_0]u\|_{L^2} \lesssim h\|\textrm{Op}^w(\{p_0,\chi_0(A_0\cdot)\}(h^{1/2}X))u\|_{L^2}+h^2\|u\|_{L^2}.
\end{equation}
From (\ref{sl7}), (\ref{sl8}) and (\ref{kk4}) we have that the principal symbol is elliptic near the supports of the two functions 
$$(1-\chi_0)(A_0\cdot),\quad  \{p_0,\chi_0(A_0\cdot)\}.$$
This yields that we may therefore get the estimate
\begin{multline}\label{eq12.3}
\|(1-M_0)u\|_{L^2}+\|\textrm{Op}^w(\{p_0,\chi_0(A_0 \cdot)\}(h^{1/2}X))u\|_{L^2} \\ =
\|(1-\chi_0)^w(A_0h^{1/2}X)u\|_{L^2}+\|\textrm{Op}^w(\{p_0,\chi_0(A_0 \cdot)\}(h^{1/2}X))u\|_{L^2} \\ \lesssim  \|Pu-hz(h)u\|_{L^2}+\mathcal{O}(h)\|u\|_{L^2}.
\end{multline}
From (\ref{eq12.4}), (\ref{bb13}) and (\ref{eq12.3}) we then have that
$$h\|u\|_{L^2} \lesssim   \|Pu-hz(h)u\|_{L^2}+  h|R_+u|,$$
when $0<h \ll 1$. Since $N_0 \geq 1$, this implies
\begin{multline}\label{eq11}
 h^{\frac{N_0}{2}+1}\|u\|_{L^2} \lesssim   h^{\frac{N_0}{2}}\|Pu-hz(h)u\|_{L^2}+  h^{\frac{N_0}{2}+1}|R_+u| \\
  \lesssim  \|Pu-hz(h)u\|_{L^2}+  h^{\frac{N_0}{2}+1}|R_+u| ,
\end{multline}
when $0<h \ll 1$. 
On the other hand, from (\ref{eq10.1}) and (\ref{eq2}) we have that  
\begin{align}\label{bb21}
& \ |E_{\pm}R_+u| \leq \Big|E_-\Big(\sum_{k=0}^{2N_0+2}a_k^w(x,D_x)h^{1+\frac{k}{2}}\Big)u\Big|+\mathcal{O}(h^{N_0+\frac{5}{2}})\|u\|_{L^2}  \\ \notag
\leq & \ |E_-(P-hz(h))u|+\mathcal{O}(h^{[\frac{N_0}{2}]+2})\|u\|_{L^2}\\ \notag
& \ \hspace{3cm}+h^{\frac{N_0}{2}+\frac{3}{2}}\sum_{|\alpha|=N_0+3}|E_-\textrm{Op}^w(X^{\alpha}R_{\alpha}(X;h))u|\\ \notag
\lesssim & \ \|(P-hz(h))u\|_{L^2}+\mathcal{O}(h^{[\frac{N_0}{2}]+2})\|u\|_{L^2}\\ \notag
&\ \hspace{3cm} +h^{\frac{N_0}{2}+\frac{3}{2}}\sum_{|\alpha|=N_0+3}|E_-\textrm{Op}^w(X^{\alpha}R_{\alpha}(X;h))u|,
\end{align}
because $\|E_-\|_{\mathscr{L}(L^2,\cc^d)}=\mathcal{O}(1)$. 
It follows from (\ref{sl72}) and (\ref{sl76}) that
\begin{align}\label{bb20}
& \ (E_-\textrm{Op}^w(X^{\alpha}R_{\alpha}(X;h))u)_{\text{\rm $k$-th component}}\\ \notag
= & \ \Big(\textrm{Op}^w(X^{\alpha}R_{\alpha}(X;h))u,\sum_{j=0}^{2N_0+2}\psi_{j,k}^-h^{\frac{j}{2}}\Big)_{L^2}\\ \notag
= & \ \sum_{j=0}^{2N_0+2}h^{\frac{j}{2}}(u,\textrm{Op}^w(X^{\alpha}\overline{R_{\alpha}(X;h)})\psi_{j,k}^-)_{L^2}=\mathcal{O}(1)\|u\|_{L^2},
\end{align}
since $\psi_{j,k}^- \in \mathscr{S}(\rr^n)$.
We therefore get from (\ref{bb21}) and (\ref{bb20}) that 
$$|E_{\pm}R_+u| \lesssim  \|Pu-hz(h)u\|_{L^2}+\mathcal{O}(h^{\frac{N_0}{2}+\frac{3}{2}})\|u\|_{L^2}+\mathcal{O}(h^{[\frac{N_0}{2}]+2})\|u\|_{L^2}.$$
If the estimate (\ref{eq10.4}) holds, we thus have
\begin{multline*}
c_0 h^{\frac{N_0}{2}+1}|R_+u| \leq |E_{\pm}R_+u| \\ \lesssim \|Pu-hz(h)u\|_{L^2}+\mathcal{O}(h^{\frac{N_0}{2}+\frac{3}{2}})\|u\|_{L^2}+\mathcal{O}(h^{[\frac{N_0}{2}]+2})\|u\|_{L^2},
\end{multline*}
and deduce from (\ref{eq11}) that
$$h^{\frac{N_0}{2}+1}\|u\|_{L^2} \lesssim \|Pu-hz(h)u\|_{L^2} +\mathcal{O}(h^{\frac{N_0}{2}+\frac{3}{2}})\|u\|_{L^2}+\mathcal{O}(h^{[\frac{N_0}{2}]+2})\|u\|_{L^2}.$$
This shows that 
$$h^{\frac{N_0}{2}+1}\|u\|_{L^2} \lesssim \|Pu-hz(h)u\|_{L^2},$$
when $0<h \ll 1$. Hence estimate (\ref{eq8}) holds true for any given Schwartz function and by density it also holds true for all $u \in L^2(\rr^n)$.
This finally proves the second implication and ends the proof of Theorem~\ref{th1}. \begin{flushright}{$\Box$}\end{flushright}

\section{Appendix}\label{appendix}
This appendix gathers miscellaneous facts and notations about quadratic differential operators used in the previous sections. We refer the reader to~\cite{kps2,kps4,kps21} as references for the results recalled in this section.

Associated with a complex-valued quadratic form
\begin{eqnarray*}
q : \rr_x^n \times \rr_{\xi}^n &\longrightarrow& \cc \\
 (x,\xi) & \mapsto & q(x,\xi),
\end{eqnarray*}
with $n \geq 1$, one has the \textit{Hamilton map} $F \in M_{2n}(\cc)$, uniquely defined by the identity
\begin{equation}
\label{10}
q((x,\xi);(y,\eta))=\sigma((x,\xi),F(y,\eta)), \quad (x,\xi) \in \rr^{2n},\  (y,\eta) \in \rr^{2n},
\end{equation}
where $q(\textrm{\textperiodcentered};\textrm{\textperiodcentered})$ stands for the polarized form
associated with the quadratic form~$q$ and where $\sigma$ is the canonical symplectic form on $\rr^{2n}$,
\begin{equation}\label{11}
\sigma((x,\xi),(y,\eta))=\xi \cdot y-x \cdot\eta, \quad (x,\xi) \in \rr^{2n},\  (y,\eta) \in \rr^{2n}.
\end{equation}
It readily follows from the definition that the real and imaginary parts of the Hamilton map
$$\textrm{Re }F=\frac{1}{2}(F+\overline{F}), \quad \textrm{Im }F=\frac{1}{2i}(F-\overline{F}),$$
$\overline{F}$ being the complex conjugate of $F$, are the Hamilton maps associated
with the quadratic forms $\textrm{Re } q$ and $\textrm{Im }q$.
The singular space $S$ associated with the quadratic symbol $q$ was introduced in~\cite{kps2} and defined as
\begin{equation}\label{h1}
S=\Big(\bigcap_{j=0}^{2n-1}\textrm{Ker}\big(\textrm{Re }F(\textrm{Im }F)^j \bigr)\Big) \bigcap \rr^{2n}.
\end{equation}
This linear subspace of the phase space plays a basic role in the understanding of the properties of the quadratic operator
$$q^w(x,D_x)u(x)=\frac{1}{(2\pi)^n}\int_{\rr^{2n}}e^{i(x-y) \cdot \xi}q\Big(\frac{x+y}{2},\xi\Big)u(y)d\xi dy,$$
when its symbol may fail to satisfy the ellipticity condition
$$(x,\xi) \in \rr^{2n}, \quad q(x,\xi)=0\, \Longrightarrow\,  (x,\xi)=0.$$ 
In particular, the known description of the spectrum of elliptic quadratic operators~\cite{S} extends to certain classes of ``partially 
elliptic'' quadratic operators. More specifically, when $q$ is a quadratic symbol with a nonnegative real part $\textrm{Re }q \geq 0$, 
satisfying the following ellipticity condition on its singular space $S$ (partial ellipticity),
\begin{equation}
\label{sm2}
(x,\xi) \in S, \quad q(x,\xi)=0\, \Longrightarrow\,  (x,\xi)=0,
\end{equation}  
then the spectrum $\mathrm{Spec}(q^w(x,D_x)$ 
of the quadratic operator  $q^w(x,D_x)$ is only composed of eigenvalues with finite algebraic multiplicities \cite{kps2} (Theorem~1.2.2) and explicitly given by
\begin{equation}\label{sm6}
\mathrm{Spec}(q^w(x,D_x))=\Big\{ \sum_{\substack{\lambda \in \mathrm{Spec}(F), \\  -i \lambda \in {\cc}_+
\cup (\Sigma(q|_S) \setminus \{0\})} }
{(r_{\lambda}+2 k_{\lambda})(-i\lambda) ; k_{\lambda} \in \nn}\Big\},
\end{equation}
where $r_{\lambda}$ is the dimension of the space of generalized eigenvectors of $F$ in $\cc^{2n}$ belonging to the eigenvalue $\lambda \in \cc$, 
and where
$$\Sigma(q|_S)=\overline{q(S)} \subset i\rr, \quad \cc_+=\{z \in \cc; \textrm{Re }z>0\}.$$
Equivalently, the singular space may be defined as the subset in phase space where all the Poisson brackets
$H_{\textrm{Im}\, q}^k\textrm{Re }q$, with $k \geq 0$, are vanishing:
$$S=\{X \in \rr^{2n} ; \ H_{\textrm{Im}\, q}^k\textrm{Re }q(X)=0,\ k \geq 0\}.$$
This shows that the singular space corresponds exactly to the set of points $X_0$ in the phase space where the real part of the symbol $q$ composed with the flow generated by the
Hamilton vector field associated with its imaginary part $\textrm{Im }q$,
$$t \longmapsto \textrm{Re }q(e^{tH_{\textrm{Im}\,q}}X_0),$$
vanishes to infinite order at $t=0$. 
Furthermore, quadratic operators with zero singular space were shown to enjoy noticeable subelliptic properties~\cite{kps21}. Namely, 
when $q$ is a complex-valued quadratic form with a nonnegative real part  $\textrm{Re }q \geq 0$, and a zero singular space  $S=\{0\}$, 
then the quadratic operator $q^w(x,D_x)$ fulfills the subelliptic estimate with a loss of $2k_0/(2k_0+1)$ derivatives
\begin{equation}\label{dl1}
\|\langle(x,D_x)\rangle^{2/(2k_0+1)} u\|_{L^2} \leq C(\|q^w(x,D_x) u\|_{L^2}+\|u\|_{L^2}), \quad u \in \mathscr{S}(\rr^n),
\end{equation}
where $\langle(x,D_x)\rangle^2=1+|x|^2+|D_x|^2$, and where $0 \leq k_0 \leq 2n-1$ stands for the smallest integer satisfying
$$\Big(\bigcap_{j=0}^{k_0}\textrm{Ker}\bigl(\textrm{Re }F(\textrm{Im }F)^j \bigr)\Big) \cap \rr^{2n}=\{0\}.$$

\bigskip

\noindent
{\bf Acknowledgements.}
The second author is most grateful for the support of the CNRS chair of excellence at Cergy-Pontoise University and the great hospitality of the University of Bologna.

\end{document}